

\def\autori{G.\ DAL MASO and F.\ MURAT}
\def\titolo{Asymptotic behaviour and correctors for Dirichlet problems}

\overfullrule=0pt

\magnification=1200
\mathsurround=1pt
\baselineskip=15pt

\topskip=15pt 
\hfuzz=2pt

\font\twelverm=cmr12
\font\twelvei=cmmi12
\font\twelvesy=cmsy10
\font\twelvebf=cmbx12
\font\twelvett=cmtt12
\font\twelveit=cmti12
\font\twelvesl=cmsl12

\font\ninerm=cmr9
\font\ninei=cmmi9
\font\ninesy=cmsy9
\font\ninebf=cmbx9
\font\ninett=cmtt9
\font\nineit=cmti9
\font\ninesl=cmsl9

\font\eightrm=cmr8
\font\eighti=cmmi8
\font\eightsy=cmsy8
\font\eightbf=cmbx8
\font\eighttt=cmtt8
\font\eightit=cmti8
\font\eightsl=cmsl8

\font\sixrm=cmr6
\font\sixi=cmmi6
\font\sixsy=cmsy6
\font\sixbf=cmbx6

\catcode`@=11 
\newskip\ttglue

\def\twelvepoint{\def\rm{\fam0\twelverm}
\textfont0=\twelverm  \scriptfont0=\ninerm
\scriptscriptfont0=\sevenrm
\textfont1=\twelvei  \scriptfont1=\ninei  \scriptscriptfont1=\seveni
\textfont2=\twelvesy  \scriptfont2=\ninesy
\scriptscriptfont2=\sevensy
\textfont3=\tenex  \scriptfont3=\tenex  \scriptscriptfont3=\tenex
\textfont\itfam=\twelveit  \def\it{\fam\itfam\twelveit}%
\textfont\slfam=\twelvesl  \def\sl{\fam\slfam\twelvesl}%
\textfont\ttfam=\twelvett  \def\tt{\fam\ttfam\twelvett}%
\textfont\bffam=\twelvebf  \scriptfont\bffam=\ninebf
\scriptscriptfont\bffam=\sevenbf  \def\bf{\fam\bffam\twelvebf}%
\tt  \ttglue=.5em plus.25em minus.15em
\normalbaselineskip=15pt
\setbox\strutbox=\hbox{\vrule height10pt depth5pt width0pt}%
\let\sc=\tenrm  \let\big=\twelvebig  \normalbaselines\rm}

\def\tenpoint{\def\rm{\fam0\tenrm}
\textfont0=\tenrm  \scriptfont0=\sevenrm  \scriptscriptfont0=\fiverm
\textfont1=\teni  \scriptfont1=\seveni  \scriptscriptfont1=\fivei
\textfont2=\tensy  \scriptfont2=\sevensy  \scriptscriptfont2=\fivesy
\textfont3=\tenex  \scriptfont3=\tenex  \scriptscriptfont3=\tenex
\textfont\itfam=\tenit  \def\it{\fam\itfam\tenit}%
\textfont\slfam=\tensl  \def\sl{\fam\slfam\tensl}%
\textfont\ttfam=\tentt  \def\tt{\fam\ttfam\tentt}%
\textfont\bffam=\tenbf  \scriptfont\bffam=\sevenbf
\scriptscriptfont\bffam=\fivebf  \def\bf{\fam\bffam\tenbf}%
\tt  \ttglue=.5em plus.25em minus.15em
\normalbaselineskip=12pt
\setbox\strutbox=\hbox{\vrule height8.5pt depth3.5pt width0pt}%
\let\sc=\eightrm  \let\big=\tenbig  \normalbaselines\rm}

\def\ninepoint{\def\rm{\fam0\ninerm}
\textfont0=\ninerm  \scriptfont0=\sixrm  \scriptscriptfont0=\fiverm
\textfont1=\ninei  \scriptfont1=\sixi  \scriptscriptfont1=\fivei
\textfont2=\ninesy  \scriptfont2=\sixsy  \scriptscriptfont2=\fivesy
\textfont3=\tenex  \scriptfont3=\tenex  \scriptscriptfont3=\tenex
\textfont\itfam=\nineit  \def\it{\fam\itfam\nineit}%
\textfont\slfam=\ninesl  \def\sl{\fam\slfam\ninesl}%
\textfont\ttfam=\ninett  \def\tt{\fam\ttfam\ninett}%
\textfont\bffam=\ninebf  \scriptfont\bffam=\sixbf
\scriptscriptfont\bffam=\fivebf  \def\bf{\fam\bffam\ninebf}%
\tt  \ttglue=.5em plus.25em minus.15em
\normalbaselineskip=11pt
\setbox\strutbox=\hbox{\vrule height8pt depth3pt width0pt}%
\let\sc=\sevenrm  \let\big=\ninebig  \normalbaselines\rm}

\def\eightpoint{\def\rm{\fam0\eightrm}
\textfont0=\eightrm  \scriptfont0=\sixrm  \scriptscriptfont0=\fiverm
\textfont1=\eighti  \scriptfont1=\sixi  \scriptscriptfont1=\fivei
\textfont2=\eightsy  \scriptfont2=\sixsy  \scriptscriptfont2=\fivesy
\textfont3=\tenex  \scriptfont3=\tenex  \scriptscriptfont3=\tenex
\textfont\itfam=\eightit  \def\it{\fam\itfam\eightit}%
\textfont\slfam=\eightsl  \def\sl{\fam\slfam\eightsl}%
\textfont\ttfam=\eighttt  \def\tt{\fam\ttfam\eighttt}%
\textfont\bffam=\eightbf  \scriptfont\bffam=\sixbf
\scriptscriptfont\bffam=\fivebf  \def\bf{\fam\bffam\eightbf}%
\tt  \ttglue=.5em plus.25em minus.15em
\normalbaselineskip=9pt
\setbox\strutbox=\hbox{\vrule height7pt depth2pt width0pt}%
\let\sc=\sixrm  \let\big=\eightbig  \normalbaselines\rm}

\def\twelvebig#1{{\hbox{$\textfont0=\twelverm\textfont2=\twelvesy
        \left#1\vbox to10pt{}\right.\n@space$}}}
\def\tenbig#1{{\hbox{$\left#1\vbox to8.5pt{}\right.\n@space$}}}
\def\ninebig#1{{\hbox{$\textfont0=\tenrm\textfont2=\tensy
        \left#1\vbox to7.25pt{}\right.\n@space$}}}
\def\eightbig#1{{\hbox{$\textfont0=\ninerm\textfont2=\ninesy
        \left#1\vbox to6.5pt{}\right.\n@space$}}}

\def\displayliness#1{\null\,\vcenter{\openup1\jot \m@th
\ialign{\strut\hfil$\displaystyle{##}$\hfil
\crcr#1\crcr}}\,}
\def\displaylinesno#1{\displ@y \tabskip=\centering
\halign to\displaywidth{ \hfil$\@lign \displaystyle{##}$ \hfil
\tabskip=\centering
&\llap{$\@lign##$}\tabskip=0pt \crcr#1\crcr}}
\def\ldisplaylinesno#1{\displ@y \tabskip=\centering
\halign to\displaywidth{ \hfil$\@lign \displaystyle{##}$\hfil
\tabskip=\centering
&\kern-\displaywidth
\rlap{$\@lign##$}\tabskip=\displaywidth \crcr#1\crcr}}

\catcode`@=12 


%

\font\sixrm=cmr6
\newcount\tagno \tagno=0                        
\newcount\thmno \thmno=0                        
\newcount\bibno \bibno=0                        
\newcount\chapno\chapno=0                       
\newcount\verno            
\newif\ifproofmode
\proofmodetrue
\newif\ifwanted
\wantedfalse
\newif\ifindexed
\indexedfalse

\def\ifundefined#1{\expandafter\ifx\csname+#1\endcsname\relax}

\def\Wanted#1{\ifundefined{#1} \wantedtrue
\immediate\write0{Wanted #1
\the\chapno.\the\thmno}\fi}

\def\Increase#1{{\global\advance#1 by 1}}

\def\Assign#1#2{\immediate
\write1{\noexpand\expandafter\noexpand\def
 \noexpand\csname+#1\endcsname{#2}}\relax
 \global\expandafter\edef\csname+#1\endcsname{#2}}

\def\pAssign#1#2{\write1{\noexpand\expandafter\noexpand\def
 \noexpand\csname+#1\endcsname{#2}}}

\def\lPut#1{\ifproofmode\llap{\hbox{\sixrm #1\ \ \ }}\fi}
\def\rPut#1{\ifproofmode$^{\hbox{\sixrm #1}}$\fi}



\def\chp#1{\global\tagno=0\global\thmno=0\Increase\chapno
\Assign{#1}
{\the\chapno}{\lPut{#1}\the\chapno}}


\def\thm#1{\Increase\thmno
\Assign{#1}{\the\chapno.\the\thmno}\the\chapno.\the\thmno
    \rPut{#1}}


\def\frm#1{\Increase\tagno
  \Assign{#1}{\the\chapno.\the\tagno}
    \lPut{#1}{\the\chapno.\the\tagno}}


\def\bib#1{\Increase\bibno
\Assign{#1}{\the\bibno}\lPut{#1}{\the\bibno}}


\def\pgp#1{\pAssign{#1/}{\the\pageno}}


\def\ix#1#2#3{\pAssign{#2}{\the\pageno}
\immediate\write#1{\noexpand\idxitem{#3}
{\noexpand\csname+#2\endcsname}}}


\def\rf#1{\Wanted{#1}\csname+#1\endcsname\relax\rPut {#1}}


\def\rfp#1{\Wanted{#1}\csname+#1/\endcsname\relax\rPut{#1}}

\input Label.aux
\Increase\verno

\immediate\write1{\noexpand\verno=\the\verno}

\ifindexed
\immediate\openout2=\jobname.idx
\immediate\openout3=\jobname.sym
\fi


\def\parag#1#2{\goodbreak\bigskip\bigskip\noindent
                   {\bf #1.\ \ #2}
                   \nobreak\bigskip}
\def\intro#1{\goodbreak\bigskip\bigskip\goodbreak\noindent
                   {\bf #1}\nobreak\bigskip\nobreak}
\long\def\th#1#2{\goodbreak\bigskip\noindent
                {\bf Theorem #1.\ \ \it #2}}
\long\def\lemma#1#2{\goodbreak\bigskip\noindent
                {\bf Lemma #1.\ \ \it #2}}
\long\def\prop#1#2{\goodbreak\bigskip\noindent
                  {\bf Proposition #1.\ \ \it #2}}
\long\def\cor#1#2{\goodbreak\bigskip\noindent
                {\bf Corollary #1.\ \ \it #2}}

\long\def\rem#1#2{\goodbreak\bigskip\noindent
                 {\bf Remark #1.\ \ \rm #2}}

\def\negbigskip{\vskip-\bigskipamount}

\def\proof{\nobreak\vskip.4cm\noindent{\it Proof.\ \ }}
\def\proofof#1{\nobreak
                \vskip.4cm\noindent{\it Proof of #1.\ \ }}

                \def\sqr#1#2{\vbox{
   \hrule height .#2pt
   \hbox{\vrule width .#2pt height #1pt \kern #1pt
      \vrule width .#2pt}
   \hrule height .#2pt }}
\def\square{\sqr74}

\def\endproof{\hphantom{MM}\hfill
\llap{$\square$}\goodbreak}

\mathchardef\emptyset="001F
\mathchardef\hyphen="002D

\def\supp{{\rm supp\,}}

\def\r{{\bf R}}
\def\rn{\r^n}
\def\rnn{\r^{n \times n}}

\def \n{{\bf N}}

\def\z{{\bf Z}}

\def\ws{w^*}


\def\e{\varepsilon}
\def\om{\omega}
\def\Om{\Omega}

\def\intK{\int_{K}}
\def\into{\int_\Omega}
\def\intu{\int_{U}}
\def\intv{\int_{V}}

\def\cp{{\rm cap}}
\def\div{{\rm div}}
\def\supp{{\rm supp}}
\def\wto{\rightharpoonup}


\def\Ae{A^{\e}}
\def\Ao{A^{0}}
\def\Abe{\overline A{}^{\e}}
\def\Abo{\overline A{}^{0}}
\def\Aie{A_{i}^{\e}}
\def\Aio{A_{i}^{0}}
\def\Aonee{A_{1}^{\e}}
\def\Aoneo{A_{1}^{0}}
\def\Atwoe{A_{2}^{\e}}
\def\Atwoo{A_{2}^{0}}
\def\Be{B^{\e}}
\def\Bie{B_{i}^{\e}}
\def\ce{c^{\e}}
\def\Ce{C^{\e}}
\def\Cie{C_{i}^{\e}}

\def\cinftyo{C^\infty_c(\Omega)}
\def\cp{{\rm cap}}

\def\deltaa{\eta}
\def\Die{D_i^{\e}}
\def\distr{{\cal D}^{\prime}(\Om)}
\def\distru{{\cal D}^{\prime}(U)}

\def\Dj#1{D_{j}#1}
\def\Eie{E_i^{\e}}
\def\ej{e_{j}}
\def\fb{\overline f}
\def\fbe{\overline f{}^{\e}}
\def\fbo{\overline f{}^{0}}
\def\fd{f_{\deltaa}}
\def\fe{f^{\e}}

\def\fo{f^{0}}
\def\fs{f}
\def\gabe{\overline\gamma{}^{\e}}
\def\gie{g^{\e}}

\def\h{H^{1}(\Om)}
\def\hc{H^{1}_{c}(\Om)}

\def\hh{H^{2}(\Om)\cap W^{1,\infty}(\Om)}
\def\hhu{H^{2}(U)\cap W^{1,\infty}(U)}
\def\ho{H^{1}_{0}(\Om)}
\def\hop{W^{1,p}_{0}(\Om)}
\def\hou{H^{1}_{0}(U)}

\def\hp{W^{1,p}(\Om)}
\def\hploc{W^{1,p}_{\rm loc}(\Om)}
\def\hm{H^{-1}(\Om)}
\def\hmp{H^{-1}(\Om)^{+}}

\def\hmq{W^{-1,q}(\Om)}

\def\hu{H^{1}(U)}

\def\lae{\lambda^{\e}}
\def\lao{\lambda^{0}}
\def\lape{\sigma^{\e}_{\scriptscriptstyle\mskip-.5\thinmuskip\oplus}}
\def\Le{L^{\e}}
\def\Lo{L^{0}}
\def\linfty{L^{\infty}(\Om)}
\def\linftymu{L^{\infty}(\Om,\mu)}

\def\linftynn{L^{\infty}(\Om,\rnn)}

\def\linftyu{L^{\infty}(U)}

\def\loneloc{L^{1}_{\rm loc}(\Om)}
\def\ltwo{L^{2}(\Om)}
\def\ltwoloc{L^{2}_{{\rm loc}}(\Om)}

\def\ltwon{L^{2}(\Om,\rn)}

\def\lmu{L^{2}(\Om,\mu)}
\def\lmue{L^{2}(\Om,\mue)}
\def\lmuoneo{L^{2}(\Om,\muoneo)}
\def\lmutwoo{L^{2}(\Om,\mutwoo)}

\def\lmuou{L^{2}(U,\muo)}

\def\lmuo{L^{2}(\Om,\muo)}
\def\lmuoloc{L^{2}_{\rm loc}(\Om,\muo)}

\def\mat{M_{\alpha}^{\beta}(\Om)}
\def\mb{{\cal M}_{\rm b}(\Om)}
\def\mm{{\cal M}(\Om)}

\def\mo{{\cal M}_{0}^{+}(\Om)}

\def\mue{\mu^{\e}} \def\muo{\mu^{0}}
\def\muoh{\hat\mu{}^{0}}
\def\muio{\mu_{i}^{0}}
\def\muoneo{\mu_{1}^{0}}
\def\mutwoo{\mu_{2}^{0}}
\def\nube{\overline\nu{}^{\e}}

\def\nue{\nu^{\e}}

\def\nuie{\nu_{i}^{\e}}
\def\nuio{\nu_{i}^{0}}
\def\nuonee{\nu_{1}^{\e}}
\def\nuoneo{\nu_{1}^{0}}
\def\nus{\nu} 

\def\nutwoo{\nu_{2}^{0}}
\def\Ome{\Omega^{\e}}
\def\ome{\omega^{\e}}
\def\omie{\omega_i^{\e}}
\def\omo{\omega^{0}}
\def\phiie{\varphi_i^{\e}}
\def\psid{\psi_{\delta}}
\def\se{\eta^{\e}}
\def\sige{\sigma^{\e}}
\def\sigo{\sigma} 
\def\sigoo{\sigma^{0}}

\def\ube{\overline u{}^{\e}}
\def\ubo{\overline u{}^{0}}
\def\ubs{\overline u}
\def\ue{u^{\e}}
\def\ued{u_{\deltaa}^{\e}}
\def\uonee{u_{\scriptscriptstyle\mskip-.5\thinmuskip\oplus}^{\e}}
\def\uoneo{u_{\scriptscriptstyle\mskip-.5\thinmuskip\oplus}^{0}}
\def\utwoe{u_{\scriptscriptstyle\mskip-.5\thinmuskip\ominus}^{\e}}
\def\uo{u^{0}}
\def\uod{u_{\deltaa}^{0}}
\def\us{u} 
\def\vbe{\overline v{}^{\e}}
\def\vbo{\overline v{}^{0}}
\def\ve{v^{\e}}
\def\ved{v_{\delta}^{\e}}
\def\vie{v_i^{\e}}
\def\vonee{v_{\scriptscriptstyle\mskip-.5\thinmuskip\oplus}^{\e}}
\def\we{w^{\e}}
\def\wbe{\overline w{}^{\e}}
\def\wbs{\overline w} 
\def\weh{\hat  w{}^{\e}}
\def\wo{w^{0}}
\def\wbo{\overline w{}^{0}}
\def\woh{\hat  w{}^{0}}
\def\wie{w_{i}^{\e}}
\def\wio{w_{i}^{0}}
\def\wonee{w_{1}^{\e}}
\def\woneo{w_{1}^{0}}
\def\ws{w} 
\def\wtwoe{w_{2}^{\e}}
\def\wtwoo{w_{2}^{0}}
\def\X{\h\cap L^{2}(\Om,\mu)}
\def\x{\ho\cap L^{2}(\Om,\mu)}
\def\xe{\ho\cap L^{2}(\Om,\mue)}
\def\Xe{\h\cap L^{2}(\Om,\mue)}
\def\xeu{\hou\cap L^{2}(U,\mue)}
\def\xie{x_i^\e}

\def\xoh{\ho\cap L^{2}(\Om,\muoh)}
\def\xo{\ho\cap L^{2}(\Om,\muo)}
\def\Xo{\h\cap L^{2}(\Om,\muo)}
\def\xou{\hou\cap L^{2}(U,\muo)}

\def\xio{\ho\cap L^{2}(\Om,\muio)}
\def\xtwoo{\ho\cap L^{2}(\Om,\mutwoo)}

\def\ye{y^{\e}}
\def\yo{y^{0}}
\def\zje{z_{j}^{\e}}
\def\zie{\zeta_{i}^{\e}}
\def\zonee{\zeta_{1}^{\e}}
\def\zoneo{\zeta_{1}^{0}}
\def\ztwoe{\zeta_{2}^{\e}}

\proofmodefalse

\def\rightheadline{\eightpoint\hfil\titolo
\hfil\tenrm\folio}
\def\leftheadline{\tenrm\folio\hfil\eightpoint
\autori \hfil}

\headline={\ifnum\pageno=0 \zeroheadline
\else\ifodd\pageno\rightheadline
\else\leftheadline\fi\fi}
\def\zeroheadline{\null\hfil}

\nopagenumbers

\pageno=0
\hsize 14truecm
\vsize 25truecm
\hoffset=0.8truecm
\voffset=-1.55truecm

\null
\vskip 2.8truecm
{\twelvepoint
\baselineskip=1.7\baselineskip
\centerline{\bf ASYMPTOTIC BEHAVIOUR AND CORRECTORS}
\centerline{\bf FOR LINEAR DIRICHLET PROBLEMS}
 \centerline{\bf WITH SIMULTANEOUSLY VARYING}
\centerline{\bf  OPERATORS AND  DOMAINS}
}
\vskip2truecm

\centerline{Gianni DAL MASO ($^{\ast}$)}
\medskip
\centerline{Fran\c cois MURAT ($^{\ast\ast}$)}
\vfil

{\eightpoint
\baselineskip=1.2\baselineskip
\centerline{\bf Abstract}
\bigskip
\noindent
We consider a sequence of Dirichlet problems in varying domains
(or, more generally, of relaxed Dirichlet problems involving measures in
$\mo$) for second order linear elliptic operators in divergence form
with varying matrices of coefficients.
When the matrices $H$-converge to a matrix $\Ao$,
we prove that there exist a subsequence
and a measure $\muo$ in $\mo$ such that the limit problem is the relaxed
Dirichlet problem corresponding to $\Ao$ and $\muo$.
We also prove a corrector result which provides an explicit
approximation of the solutions in the $H^{1}$-norm, and which is
obtained by multiplying the corrector for the $H$-converging matrices
by  some special test function which depends both on
the varying matrices and on 
the varying domains.

\vfil
\item{($^\ast$)}SISSA, via Beirut 4, 34014 Trieste, Italy
\par
\item{}e-mail: {\tt dalmaso@sissa.it}
\par
\smallskip
\item{($^{\ast\ast}$)}Laboratoire Jacques-Louis Lions,
Bo\^\i te courrier 187, Universit\'e Pierre et Marie Curie,
75252 Paris cedex 05, France
\par
\item{}e-mail: {\tt murat@ann.jussieu.fr}

}
\vfil
\vskip 1truecm
\centerline {Ref. S.I.S.S.A.
40/02/M (May 2002)}
\vskip 1truecm
\eject
\null
\pageno=0
\eject
%
%
%
\pageno=1
\topskip=25pt
\vsize 22.5truecm
\hsize 16.2truecm
\hoffset=0truecm
\voffset=0.5truecm
\baselineskip=15pt
\hfuzz=2pt
\parindent=2em
\mathsurround=1pt

\parag{\chp{intro}}{Introduction}

In this paper we consider a sequence of
linear Dirichlet
problems
$$
\cases{\ue\in H^{1}_{0}(\Ome)\,,&
\cr\cr
-\div(\Ae D\ue ) = f \quad \hbox{in }\, {\cal D}'(\Ome) \,,&
\cr}\leqno(\frm{0.0})
$$
where the matrices $A^\varepsilon$ and the domains
$\Omega^\varepsilon$ both depend on the parameter $\e$.
We assume that the open sets
$\Omega^\varepsilon$ are all
contained in a fixed bounded open subset $\Omega$ of $\rn$, and
that the matrices $A^\varepsilon$, defined on
$\Omega$ with measurable coefficients, are
coercive and bounded, uniformly  with respect to~$\varepsilon$. Our goal is
to study the behaviour of the solutions $u^\varepsilon$ as
$\varepsilon$ tends to zero.

In the special case $\Ome=\Om$ it is known (see Section~\rf{hcon}) that there
exist  a subsequence, still
denoted by $(A^\varepsilon)$, and a matrix $\Ao$,
called the $H$-limit of
$(A^\varepsilon)$,  such that  for every
$f
\in H^{-1}(\Omega)$ the solutions $v^\varepsilon$ of
the problems
$$
\cases{\ve \in \ho \,,&
\cr\cr
\displaystyle
-\div(\Ae D\ve)= f \quad\hbox{in }\distr \,,&
\cr}
$$
converge weakly in $H^1_0(\Omega)$ to the solution
$v^0$ of
$$
\cases{v^0 \in \ho \,,&
\cr\cr
\displaystyle
-\div(A^0 D v^0)= f \quad\hbox{in }\distr \,,&
\cr}
$$
and satisfy also
$$
\Ae D\ve \wto \Ao D v^0 \qquad\hbox{weakly in }\,\ltwon\,.
$$

Without making any further hypothesis on the
open sets $\Omega^\varepsilon$, we prove in the
present paper that there exists a subsequence, still
denoted by $(\Omega^\varepsilon)$, such that for every $f
\in H^{-1}(\Omega)$ the solutions $u^\varepsilon$ of
(\rf{0.0}) converge to the solution $u^0$ of the problem
$$
\cases{\uo\in\xo\,,&
\cr\cr
\displaystyle
\into\Ao D\uo Dy\,dx +
\into \uo y \,d\muo
=\langle f, y \rangle \qquad \forall y\in\xo\,,&
\cr}\leqno(\frm{0.1})
$$
where $\mu^0$ belongs to $\mo$, a class of a nonnegative Borel measures which
vanish on all sets of capacity zero, but can take the value
$+\infty$ on some subsets of $\Omega$ (see Section~\rf{prel}).

Problems like (\rf{0.1}) are called relaxed Dirichlet problems
(see Section~\rf{relax}) and have been extensively studied to
describe the limits of the solutions of
(\rf{0.0}) when the matrices $A^\varepsilon$
do not depend on~$\e$.
On the other hand, problems (\rf{0.0}) can be written as
relaxed Dirichlet problems (see Remark \rf{bc01}) by considering the measures
$\mu^\varepsilon$ defined by
$$
\mue(B)=\cases{0,\qquad &if $\cp(B\setminus \Ome)=0$,
\cr\cr
+\infty,&otherwise.\cr}\leqno(\frm{bbbb1})
$$
Actually in the paper
we consider not only the
case of Dirichlet problems (\rf{0.0}), which correspond to the measures $\mue$
defined by (\rf{bbbb1}), but more in general we study
the case of a sequence of
relaxed Dirichlet problems with arbitrary $\mue\in\mo$.

In the limit problem (\rf{0.1})
the measure $\mu^0$ does not depend on $f$, but, as shown
in Section~\rf{examples}, it depends both on the sequence of
sets
$(\Omega^\varepsilon)$ and on the sequence of matrices
$(A^\varepsilon)$ (and not only on its $H$-limit $A^0$).
Nevertheless the sequence $(\Omega^\varepsilon)$ has a
stronger influence than the sequence $(A^\varepsilon)$:
indeed the limit measures corresponding to the same
sequence $(\Omega^\varepsilon)$ but to different
sequences $(\Aie)$ are equivalent (see Theorem~\rf{aa04}).

In Section~\rf{ho} we give a fairly general and
flexible method to construct the limit measure $\mu^0$
using suitable test functions $\omega^\varepsilon$
associated to $\Omega^\varepsilon$ and $A^\varepsilon$.
We then pass to the limit in the sequence of problems
(\rf{0.0}) by a duality argument and obtain (\rf{0.1}).

In Section~\rf{corrector} we continue the study of the
behaviour of the solutions
$u^\varepsilon$ of (\rf{0.0}) by giving a corrector result.
By this we mean the following: when the solution
$u^0$ of the limit problem (\rf{0.1}) can be written as
$$
u^0 = \psi\, \omega^0\,,\leqno(\frm{0.3})
$$
where $\omega^0$ is the limit of the above test
functions $\omega^\varepsilon$ and $\psi$ is
sufficiently smooth (actually in $H^2(\Omega) \cap
W^{1,\infty}(\Omega)$), we prove that
$$
\ue=(\psi + \sum_{j=1}^{n} \Dj{\psi}\, \zje)\,
\ome +r^\varepsilon
\qquad\hbox{with }r^\varepsilon \rightarrow 0\,\hbox{ strongly in
}\,H^1_0(\Omega)\,,
\leqno(\frm{af1})
$$
where the functions $z^\varepsilon_j$
depend only on the matrices~$A^\varepsilon$.
This provides an approximation of
$u^\varepsilon$ in the norm of
$H^1_0(\Omega)$ by means of functions that are
constructed explicitly.

When (\rf{0.3})
is not satisfied with a smooth $\psi$, a similar but more
technical result holds (see Theorem~\rf{ab03}). We
also prove a local version of this corrector result.

Moreover, we prove (global and local) convergence
and corrector results when also the right hand side
of (\rf{0.0}) depends on $\e$
and converges strongly in a
convenient sense (see Section~\rf{Le}).

Let us finally note that results similar to those presented in this
paper have been recently obtained by Calvo Jurado and
Casado Diaz in [\rf{Cal-Cas}] for a class of nonlinear monotone
elliptic equations.

\goodbreak\bigskip
\noindent{\it Contents\/}\nobreak
\item{1.}Introduction

\item{2.}Preliminaries on capacity and measures

\item{3.}$H$-convergence

\item{4.}Relaxed Dirichlet problems

\item{5.}A global convergence result

\item{6.}An example

\item{7.}Global and local corrector results

\item{8.}A comparison theorem

\item{9.}Proofs of the corrector results

\item{10.}Problems with more general data

\parag{\chp{prel}}{Preliminaries on capacity and measures}

In this section we first introduce a few notation.
Then we recall some known results
on measures,  capacity, and fine properties of Sobolev
functions.

\goodbreak\bigskip
\noindent{\it Notation\/}\nobreak

Throughout  the paper $\Omega$ is a bounded open subset of $\rn$,
$n\ge 1$. The space $\distr$ of distributions in~$\Om$ is the dual of
the space $\cinftyo$. The space $\hop$, $1\le p <+\infty$, is the
closure of $\cinftyo$ in the Sobolev space $\hp$, while
$W^{-1,q}(\Om)$, $1\le q<+\infty$, is the space of all distributions of
the
form $f=f_{0}+\sum_{j} D_{j} f_{j}$, with $f_{0}, f_{1}, \ldots,
f_{n}\in L^{q}(\Om)$ (if $1/p+1/q=1$, then $\hmq$ is the dual of
$\hop$).
In the Hilbert
case $p=q=2$ these spaces are denoted by $\ho$, $\h$, and $\hm$,
respectively. The norm in $\ho$ is defined by
$$
\|u\|_{\ho}=\Big( \into |Du|^{2}dx \Big)^{1\over 2}\,,
$$
while the duality pairing between $\hm$ and $\ho$ is denoted by
$\langle\cdot,\cdot\rangle$.
We shall sometimes use also the Sobolev space $H^{2}(\Om)=W^{2,2}(\Om)$.

The {\it adjoint\/} of a matrix $A$ is denoted by ${\overline A}$.
Since complex numbers are not used in this paper, {\it the bar never
denotes complex conjugation\/}. If $w$ is an object related to the
matrix $A$, then ${\overline w}$ denotes the corresponding object
related to the adjoint ${\overline A}$.

Throughout the paper $\e$ varies in a stricly decreasing sequence of
positive real numbers which converges to $0$.
When we write $\e>0$, we consider only the elements of this sequence,
while when we write $\e\ge0$ we also consider its limit $\e=0$.

\goodbreak\bigskip
\noindent{\it Capacity and measures\/}\nobreak

For every subset $E$ of $\Om$ the {\it capacity\/} of  $E$
in $\Om$, denoted by $\cp(E)$, is defined as the infimum of
$\into |Du|^2\,dx$ over the set of all functions $u\in\ho$ such that
$u\geq 1$ a.e.\ in a neighbourhood of~$E$.
We say that a property ${\cal P}(x)$ holds {\it quasi everywhere}
(abbreviated as {\it q.e.\/}) in a set $E$ if it holds for all $x\in E$
except for a
subset $N$ of $E$ with  ${\cp(N)=0}$. The expression {\it almost
everywhere\/} (abbreviated as {\it a.e.\/}) refers, as usual, to the
analogous
property for the Lebesgue measure.

A function $u\colon\Om\to\r$ is said to be {\it quasi continuous\/} if
for every $\e>0$ there exists a set $E\subseteq\Om$, with
$\cp(E)<\e$, such that the restriction of $u$ to
${\Om\setminus E}$ is continuous. A subset $U$ of $\Om$ is said to
be  {\it quasi open\/} if for every $\e>0$ there exists an open set
$V\subseteq\Om$, with $\cp(V{\scriptstyle\triangle} \, U)<\e$,
where ${\scriptstyle\triangle}$ denotes the symmetric difference.

Every $u\in\h$ has a
{\it quasi continuous representative\/}, which is uniquely defined up to
a set
of capacity zero. In the sequel we shall always identify $u$ with its
quasi
continuous representative, so that the pointwise values of a function
$u\in\h$ are  defined quasi everywhere in $\Om$.
If $u\in\h$, then
$$
u\ge0\quad \hbox{a.e.\ in }\,\Om\quad
\Longleftrightarrow
\quad
u\ge0\quad \hbox{q.e.\ in }\,\Om\,.
\leqno(\frm{bg60})
$$
If a sequence $(u_j)$  converges to $u$ strongly in $\ho$,
then a subsequence of $(u_j)$ converges to $u$ q.e.\ in $\Om$.
For all these  properties concerning quasi continuous
representatives of Sobolev functions we refer to~[\rf{Eva-Gar}],
Section~4.8,   [\rf{Maz}], Section 7.2.4, [\rf{Hei-Kil-Mar}],
Section~4,
or~[\rf{Zie}], Chapter~3.

The characteristic function $1_{E}$ of a set $E\subseteq\Om$ is defined
by $1_{E}(x)=1$ if $x\in E$ and $1_{E}(x)=0$ if $x\in{\Om\setminus E}$.
The following lemma (see [\rf{DM-83}], Lemma~1.5, or
[\rf{DM-Gar-94}], Lemma~1.1) concerns
the pointwise approximation of the
characteristic function of a quasi open set.

\lemma{\thm{be1.1}}{For every quasi open set $U$ of~$\Om$ there
exists an increasing sequence $(z_{k})$ of
nonnegative functions of $\ho$ converging to $1^{}_U$ pointwise
q.e.\ in~$\Om$.
}

\bigskip

By a {\it  nonnegative Borel measure\/} on
$\Om$ we mean a countably  additive set function defined on the Borel
subsets of $\Om$ with values in $[0,+\infty]$.
By a {\it  nonnegative Radon measure\/} on $\Om$ we  mean a
nonnegative Borel measure which is finite on every compact subset
of $\Om$. Every nonnegative Borel measure $\mu$ on $\Omega$ can be
extended to a Borel regular outer measure on $\Om$ by setting for every
subset $E$ of $\Om$
$$
\mu(E)=\inf\{\mu(B):B\,\hbox{ Borel, }\, E\subseteq B \subseteq \Om\}\,.
$$
If $\mu$ is a nonnegative Borel measure
on~$\Om$, we
shall use
$L^r(\Om,\mu)$, ${1\le r\le{+}\infty}$, to denote the usual Lebesgue
space with respect to the measure~$\mu$. We adopt the
standard notation $L^r(\Om)$ when $\mu$ is the Lebesgue measure.

We will consider the cone $\mo$ of all nonnegative Borel measures $\mu$
on $\Om$ such  that
\smallskip\item{(a)}
$\mu(B)=0$ for every Borel set $B\subseteq \Om$ with
$\cp(B)=0$,
\smallskip\item{(b)}
$\mu(B)=\inf\{\mu(U):U\ \hbox{quasi open}\,,\>B\subseteq U\}$
for every Borel set $B\subseteq \Om$.
\smallskip\noindent
If $E\subseteq \Om$ and $\cp(E)=0$, then $E$ is contained in a
Borel set $B\subseteq \Om$ with
$\cp(B)=0$. Therefore $E$ is $\mu$-measurable by (a).
Property (b) is a weak regularity property of the measure~$\mu$.
It is always satisfied if $\mu$ is a nonnegative Radon measure.
Since any quasi open set differs from a Borel set by a set of
capacity zero, every quasi open set is
$\mu$-measurable for every nonnegative Borel measure~$\mu$ which
satisfies~(a).

Let us explicitly observe that the notation is not fixed in the
literature and that in other works (see, e.g.,~[\rf{DM-Mos-87}])
${\cal M}_0(\Om)$
denotes the set of nonnegative Borel measures
which only satisfy~(a), while the
set that we call $\mo$ in the present paper is sometimes denoted by
${\cal M}^*_0(\Om)$ (see, e.g.,~[\rf{DM-87}]).

For every quasi open set $U\subseteq\Om$ we define the Borel
measure $\mu_U$ by
$$
\mu_U(B)=\cases{0,\qquad &if $\cp(B\setminus U)=0$,
\cr\cr
+\infty,&otherwise.\cr}\leqno(\frm{bb1})
$$
Roughly speaking, $\mu_U$ is identically zero on $U$ and
identically $+\infty$ on ${\Om\setminus U}$.
It is easy to see that this measure belongs to the
class $\mo$. Indeed, property (a) follows immediately from the
definition, and it is enough to verify (b) only for every Borel set
with $\mu_U(B)<+\infty$; in this case $\cp(B\setminus U)=0$,
and this implies that $V=U\cup B$ is quasi open (since $U$ is quasi
open),
contains $B$, and $\mu_U(V)=0$ (since
$\cp(V\setminus U)=\cp(B\setminus U)=0$),
so that (b) is satisfied. The measures $\mu_U$ will be used to
transform a sequence of Dirichlet problems on varying domains into a
sequence of relaxed Dirichlet problems on a fixed domain
(see Remark~\rf{bc01} and the proof of Corollary~\rf{bi01}).

If $\mu\in\mo$, then the space $\h\cap L^{2}(\Om,\mu)$ is well defined,
since every function $u$ in $\h$ is defined $\mu$-almost everywhere
and is $\mu$-measurable in~$\Om$ (recall that $u$ is quasi continuous,
so
that $\{u>t\}$ is quasi open for every $t\in\r$). It
is easy to see that $\h\cap L^{2}(\Om,\mu)$ is a Hilbert space for the
scalar product
$$
(u,v)_{\h\cap L^{2}(\Om,\mu)}= \into Du\,Dv\,dx +
\into u\,v\,dx + \into u\,v\,d\mu
\leqno(\frm{bd5})
$$
(see [\rf{But-DM-91}], Proposition~2.1).

The space of all (signed) Radon measures on $\Om$ will be denoted by
$\mm$,
while $\mb$ will be the space of all $\mu\in\mm$ with
$|\mu|(\Om)<+\infty$,
where $|\mu|$ denotes the total variation of $\mu$.
A subset $\cal A$ of $\mm$ is bounded if for every compact set
$K\subseteq\Om$ we have
$$
\sup_{\mu\in{\cal A}} |\mu|(K) <+\infty\,.
$$

Every Radon measure
on $\Om$ will be identified with an element of
$\distr$ in the usual way. Therefore
$\mu$ belongs to $\mm \cap W^{-1,q}(\Om)$ if and only if there exist
$f_{0}, f_{1}, \ldots, f_{n}\in L^{q}(\Om)$ such that
$$
\into \varphi \,d\mu = \into f_{0} \varphi \,dx
- \sum_{j=1}^{n} \into f_{j} D_{j} \varphi \,dx \qquad
\forall \varphi\in\cinftyo\,.
$$
Note that, by the Riesz  theorem, every nonnegative element of
$W^{-1,p}(\Om)$ is a nonnegative  Radon measure  on $\Om$.

The cone
of all nonnegative elements of $\hm$ will be denoted by $\hmp$.
It is well known that every element of $\hmp$ is a
nonnegative Radon measure which belongs also to $\mo$. In other
words we have the inclusion $\hmp\subseteq \mm\cap\mo$.

\parag{\chp{hcon}}{$H$-convergence}

In this section we recall the definition of
$H$-convergence and the
corresponding corrector result. Moreover we prove a
fairly general convergence
theorem for right hand sides which do not converge strongly in
$\hm$.

Throughout the paper we fix two constants
$\alpha$ and $\beta$ such that
$$
0 < \alpha \leq \beta < + \infty\,.
$$
We define $\mat$ as the set of all matrices $A$ in $\linftynn$
such that
$$
A(x) \ge \alpha I\,, \qquad (A(x))^{-1} \ge \beta^{-1} I\,,
\qquad\hbox{for a.e.\ } \, x \in \Om\,.
\leqno(\frm{za12})
$$
In (\rf{za12}) $I$ is the identity matrix in $\rnn$,
and the inequalities are in the sense of the quadratic forms
defined by$A(x) \xi \,\xi$ for $\xi \in \rn$.
Note that (\rf{za12}) implies that
$$
|A(x)| \le \beta  \qquad\hbox{for a.e.\ }\, x \in \Om\,,
\leqno(\frm{za13})
$$
and that necessarily $\alpha \le \beta$.

\goodbreak\bigskip
\noindent{\it Definition of $H$-convergence}\nobreak

A sequence $(\Ae)$ of matrices in $\mat$
{\it $H$-converges\/} to a matrix $\Ao$ in $\mat$
if for every $f \in \hm$ the sequence $(\ue)$ of the solutions to the
problems
$$
\cases{\ue \in \ho \,,&
\cr\cr
\displaystyle
-\div(\Ae D\ue)= f \quad\hbox{in }\distr \,,&
\cr}\leqno(\frm{za14})
$$
satisfies
$$
\displayliness{
\ue \wto \uo \qquad\hbox{weakly in }\, \ho \,,
\cr
\Ae D\ue \wto \Ao D\uo \qquad\hbox{weakly in }\,\ltwon\,,
\cr}
$$
where $\uo$ is the solution to the problem
$$
\cases{\uo \in \ho \,,&
\cr\cr
\displaystyle
-\div(\Ao D\uo)= f \quad\hbox{in }\distr \,.&
\cr}\leqno(\frm{za16})
$$

Every sequence of matrices in $\mat$
has a subsequence
which $H$-converges to a matrix in $\mat$ (see [\rf{Spa}] and
[\rf{Mur-Tar}]).


Denoting the adjoint of $\Ae$ by $\Abe$,
it is easy to prove that the sequence $(\Abe)$ $H$-converges to $\Abo$
when the sequence $(\Ae)$ $H$-converges to $\Ao$.

 If $U$ is an open set contained in $\Om$,
we can consider also the
notion of  $H$-convergence in $U$,
replacing $\Om$ by $U$ in the
definition. It is not difficult to prove that
$(\Ae)$ $H$-converges to $\Ao$
in $U$, for every open set
$U\subseteq\Om$, if
$(\Ae)$ $H$-converges to $\Ao$ in $\Om$.

%
%

\goodbreak\bigskip
\noindent{\it Corrector result}\nobreak

Besides the compactness result mentioned above,
one of the main theorems is the corrector result
(see [\rf{Mur-Tar}], and [\rf{BLP}], [\rf{Sanchez}] in the periodic
case).
Let $(e_{1}, e_{2},  \ldots, e_{n})$ be the canonical basis of $\rn$.
For $j = 1, 2, \ldots, n$ there exists a sequence
$(\zje)$ in $\h$ such that
$$
\ldisplaylinesno{
\zje\wto 0 \qquad\hbox{weakly in }\, \h\,,
&(\frm{ab3})
\cr
\Ae(D\zje+\ej) \wto \Ao\ej  \qquad\hbox{weakly in }\, \ltwon \,,
&(\frm{ab6})
\cr
-\div(\Ae(D\zje+\ej)) \to -\div(\Ao\ej) \qquad\hbox{strongly in }\,
\hm\,.
&(\frm{ab2})
\cr}
$$
Throughout the paper we will also assume that
$$
\ldisplaylinesno{
\zje\to 0 \qquad\hbox{strongly in }\, \linfty\,,
&(\frm{ab4})
\cr
\zje\wto 0 \qquad\hbox{weakly in }\, \hp\, \hbox{ for some }\, p>2 \,;
&(\frm{ab5})
\cr}
$$
using De Giorgi's and Meyers' regularity theorems,
such a sequence can be constructed, for instance, by solving the
problems
$$
\cases{\zje\in H^{1}_{0}(\Om')\,,&
\cr\cr
-\div(\Ae(D\zje+\ej))= -\div(\Ao\ej) \quad\hbox{in }\, {\cal
D}'(\Om')\,,
&\cr}
$$
where $\Om'$ is a bounded open set with $\Om \subset \subset \Om'$,
and $\Ae$ is extended by $\alpha I$ on ${\Om'\setminus\Om}$.
(The use of $\Om'$ is needed here only to obtain
a global $\hp$ bound for $\zje$ in the case where $\partial \Om$
is not smooth.)

Let $f\in\hm$, let $(\ue)$ be the sequence of the solutions to
(\rf{za14}), and let $\uo$  be the solution to (\rf{za16}). Given
$\delta>0$, let $\psid$ be a function in $\cinftyo$ which satisfies
$$
\beta \into |D\uo-D\psid|^{2} dx  <  \delta \,,
\leqno(\frm{za21})
$$
and let $\ved$ be defined  by
$$
\ved=\psid + \sum_{j=1}^{n} \Dj{\psid}\,\zje \,.
\leqno(\frm{za22})
$$
Then (see [\rf{Mur-Tar}])
$$
\limsup_{\e\to0} \;\alpha \into  |D\ue-D\ved|^{2} dx
< \delta\,.
\leqno(\frm{za23})
$$

\bigskip

If $\uo$ belongs to $\cinftyo$, we can take $\psid = \uo$ in (\rf{za21})
for every $\delta > 0$, so that
$$
\ved = \ve=\uo + \sum_{j=1}^{n} \Dj{\uo}\zje \,,
\leqno(\frm{za24})
$$
and  (\rf{za23})  implies that
$$
D\ue - D\ve \to 0 \qquad\hbox{strongly in }\,\ltwon \,,
\leqno(\frm{za25})
$$
which means that $D\ue$ is equivalent to $D\ve$
(and also to
$D\uo + \sum_j \Dj{\uo} \, D{\zje} $,  using (\rf{ab3})),
as far as convergences in $\ltwon$ are concerned.

In the general case where $\uo$ only belongs to $\ho$,
we obtain from (\rf{za23}) that
$$
D\ue = D\psid + \sum_{j=1}^{n} \Dj{\psid} \, D{\zje}
+R^{\e}_{\delta}\,,
\qquad\hbox{with}\qquad
\limsup_{\e\to0}\|R^{\e}_{\delta}\|_{\ltwon}^{2} <
{\delta\over\alpha}\,.
$$
This is a corrector result:
indeed it allows one to replace $D\ue$ by an explicit expression,
up to a remainder $R^{\e}_{\delta}$ which is small in $\ltwon$
for $\delta$ small, uniformly in $\e$;
similar corrector results have been obtained also in the case of
local solutions. Applications can be found,
e.g., in [\rf{BBDM}] and [\rf{Aglio-Mur}].

\goodbreak\bigskip
\noindent{\it A convergence result}\nobreak

We conclude this section with the following convergence result,
which is implicitly used in various works
(see, e.g., [\rf{Boc-Mur 1}]). Observe that there is no boundary
condition on the solutions $\ue$ and that the right hand sides $\fe$
do not converge strongly in $\hm$.

\th{\thm{za006}}{Let $(\Ae)$ be a sequence of matrices in $\mat$
which $H$-converges to a matrix $\Ao$ in $\mat$,
and let $(\ue)$  be a sequence in $\h$ such that
$$
\displayliness{
\ue \wto \uo \qquad\hbox{weakly in }\, \h\,,
\cr
-\div(\Ae D\ue)=\fe \quad\hbox{in }\, \distr\,
\hbox{ for every }\,\e\ge 0\,.
\cr}\leqno(\frm{ae1})
$$
Assume that $\fe = \gie +\mue+\nue$ for every $\e>0$,
where
$$
\displayliness{
(\gie) \quad\hbox{is relatively compact in }\, W^{-1,p}_{\rm loc}(\Om)\,
\hbox{ for some }\,p>1\,,
\cr
(\mue) \quad\hbox{is bounded in }\, \mm\,,
\cr
\nue\ge 0 \quad\hbox{in }\, \distr \,.
\cr}\leqno(\frm{ae2})
$$
Then
$$
\displayliness{
\fe \wto \fo \, \hbox{ weakly in }\, \hm \,
\hbox{ and strongly in }\, W^{-1,q}_{\rm loc}(\Om)
\, \hbox{ for every } \, q<2\,,
\cr
\Ae D\ue \wto \Ao D\uo \qquad\hbox{weakly in }\,\ltwon\,.
\cr}\leqno(\frm{ae3})
$$
}

In the present paper, this theorem will be used
with $\mue = 0$ and $(\gie)$ relatively compact (or even constant) in
$\hm$.

\proof
Let $K$  be any compact set of $\rn$ with $K \subseteq \Om$,
and let $\varphi\in\cinftyo$ with $\varphi \ge 0$ on $\Om$
and $\varphi = 1$ on $K$. We have
$$
0 \le \intK d\nue \le \into\varphi \, d\nue
= \into\Ae D\ue D\varphi \,dx -
\langle \gie, \varphi \rangle  - \into\varphi \, d\mue \,.
\leqno(\frm{za31})
$$
Because of (\rf{za13}), (\rf{ae1}), and (\rf{ae2}),
the right hand side of (\rf{za31}) is bounded independently of
$\e$.
This implies that
$$
(\nue) \quad\hbox{is bounded in }\, \mm\,.
\leqno(\frm{za32})
$$

For every bounded open set $U$ of $\rn$,
the embedding
$W^{1,r}_{0}(U) \subseteq C^{0}_{0}(U)$
is compact for every $r > n$.
This implies that the embedding
${\cal M}_{b}(U)  \subseteq W^{-1,s}(U)$
is compact for every $s < n/ (n - 1)$,
and therefore the embedding
$\mm \subseteq W^{-1,s}_{\rm loc}(\Om)$
is compact for every $s < n/ (n - 1)\,.$
Therefore (\rf{ae2}) and (\rf{za32}) imply that
$(\mue + \nue)$ is relatively compact in
$W^{-1,s}_{\rm loc}(\Om)$ ,
which implies that $(\fe)$ is relatively compact in $W^{-1,t}_{\rm loc}(\Om)$
for some $t>1$.
On the other hand, we deduce from
(\rf{ae1}) and (\rf{za13}) that $(\fe)$ is bounded in $\hm$.
By interpolation, $(\fe)$ is relatively compact in
$W^{-1,q}_{\rm loc}(\Om)$
for every $q<2$.

Let now $\vbo$ be an arbitrary function in $\cinftyo$,
and, for every $\e>0$, let $\vbe$ be the solution to the problem
$$
\cases{\vbe \in \ho \,,&
\cr\cr
\displaystyle
-\div(\Abe D\vbe)= -\div(\Abo D\vbo)  \quad\hbox{in }\distr \,.
&\cr}\leqno(\frm{za35})
$$
Recall that the sequence $(\Abe)$ $H$-converges to $\Abo$,
so that
$$
\displayliness{
\vbe \wto \vbo \qquad\hbox{weakly in }\, \ho\,,
\cr
\Abe D\vbe \wto \Abo D\vbo \qquad\hbox{weakly in }\, \ltwon\,,
\cr
\vbe \wto \vbo \qquad\hbox{weakly in }\, \hploc \,
\hbox{ for some }\,p > 2 \,,
\cr}\leqno(\frm{za36})
$$
where in the last assertion we have used Meyers' regularity result
(see~[\rf{Mey}]).

Let $\varphi \in \cinftyo$.
Using $\vbe\varphi$ as test function in (\rf{ae1}),
and $\ue\varphi $ as test function in (\rf{za35}),
we have
$$
\displayliness{
\langle \fe, \vbe \varphi\rangle =
\into\Ae D\ue D\vbe \varphi\,dx +
\into\Ae D\ue D\varphi \, \vbe dx=
\cr
= \langle -\div(\Abo D\vbo), \ue \varphi \rangle
- \into\Abe D\vbe D\varphi \, \ue \,dx
+ \into\Ae D\ue D\varphi \, \vbe dx  \,.
\cr } \leqno(\frm{za37})
$$
Passing to a subsequence, we may assume that
$$
\displayliness{
\Ae D\ue \wto \sigo  \qquad\hbox{weakly in }\, \ltwon\,,
\cr
\fe \wto \fs \qquad\hbox{weakly in }\,\hm\,
\hbox{ and strongly in }\, W^{-1,q}_{\rm loc}(\Om) \,
\hbox{ for every }\, q<2\,,
\cr} \leqno(\frm{za38})
$$
for some $\sigo\in\ltwon$ and $\fs\in\hm$.
It is now  easy to pass to the limit in
the left and right hand sides of
(\rf{za37}) by using (\rf{za36}), (\rf{za38}),
and Rellich's compactness theorem.
One obtains
$$
\eqalign{
\langle \fs, &  \vbo \varphi \rangle =
\cr
&= \langle -\div(\Abo D\vbo),   \uo \varphi\rangle
- \into\Abo D\vbo D\varphi \,\uo dx
+ \into \sigo\, D\varphi \, \vbo dx
\cr
&= \into\Abo D\vbo D\uo  \varphi \,dx
+ \into \sigo\, D\varphi \, \vbo dx
\cr
&= \into\Ao D\uo D\vbo \varphi \,dx
+ \into \sigo\, D\varphi \, \vbo dx  \,.
\cr } \leqno(\frm{za39})
$$
Since
$$
-\div(\sigo)=\fs \qquad\hbox{in }\, \distr\,,
\leqno(\frm{bi1})
$$
one deduces from (\rf{za39}) that
$$
\into \sigo\, D \vbo \varphi \,dx
= \into\Ao D\uo D\vbo \varphi \,dx \,,
\leqno(\frm{za40})
$$
for every $\varphi \in \cinftyo$ and every $\vbo \in \cinftyo$.
Since, for every point $x \in \Om$,
the vector $D\vbo (x)$ can be chosen
to coincide with any prescribed vector of $\rn$,
(\rf{za40}) implies that
$$
\sigo = \Ao D \uo \,\hbox{ a.e.\ in }\, \Om \,,
$$
which, together with (\rf{bi1}), gives $\fs=\fo$.
The uniqueness of the limits in (\rf{za38})
implies that the whole sequences converge,
and this completes the proof of (\rf{ae3}).
\endproof

\parag{\chp{relax}}{Relaxed Dirichlet problems}

In this section we recall the definition, introduced in [\rf{DM-Mos-86}]
and
[\rf{DM-Mos-87}], of relaxed Dirichlet problems associated with measures
$\mu\in\mo$, and prove that, under some conditions on the data,
the measure
$\mu$ can be reconstructed from the solution of the corresponding relaxed
Dirichlet problem.

\goodbreak\bigskip
\noindent{\it Relaxed Dirichlet problems}\nobreak

Given $A\in\mat$, $\mu\in\mo$, and $f\in\hm$,
we call {\it relaxed Dirichlet problem\/}
the problem of finding $u$ such that
$$
\cases{u\in\x\,,&
\cr\cr
\displaystyle
\into A\, Du\, Dy\,dx +
\into u\, y \,d\mu
=\langle f, y \rangle \qquad \forall y\in\x\,.&
\cr}\leqno(\frm{bc1})
$$

By a straightforward application of the Lax-Milgram
lemma problem (\rf{bc1}) has a unique
solution $u$ (see [\rf{DM-Mos-87}], Theorem~2.4) and $u$
satisfies the estimate
$$
\alpha \into | Du |^{2} dx + \into |u |^{2} d\mu \le
{1\over \alpha} \|f\|_{\hm}^{2} \,.
\leqno(\frm{bc10})
$$

A connection between classical Dirichlet problems on open subsets
of~$\Om$ and relaxed  Dirichlet problems of the form (\rf{bc1}) is
given by the following remark.

\rem{\thm{bc01}}{Using Theorem~4.5 of [\rf{Hei-Kil-Mar}] it is easy
to check that,
if $U\subseteq\Om$ is open and $\mu_{U}$ is the measure
introduced in~(\rf{bb1}), then $u\in\ho\cap L^{2}(\Om,\mu_{U})$
if and only if the restriction
of $u$ to $U$ belongs to $\hou$ and $u=0$ q.e.\ in ${\Om\setminus U}$.
Therefore when $\mu=\mu_{U}$ problem (\rf{bc1}) reduces to the
following boundary value
problem on $U$:
$$
\cases{u\in \hou\,,&
\cr
\cr
\displaystyle-\div(A\, D u)= f \quad\hbox{in }\distru \,,&
\cr}\leqno(\frm{bc2})
$$
in the sense that $u$ is the solution of (\rf{bc1}) if and only if its
restriction to $U$ is the solution of (\rf{bc2})
and $u=0$ q.e.\ in ${\Om\setminus U}$.

The name ``relaxed Dirichlet problem'' is motivated by the fact that
the limit of the solutions to Dirichlet problems on varying domains
$\Ome$ always satisfies a relaxed Dirichlet problem (see, e.g.,
[\rf{DM-Mos-87}] and [\rf{DM-Gar-94}], and also Corollary \rf{bi01}
below).
Moreover, the results proved in [\rf{DM-Mos-87}] and
[\rf{DM-Mal}] ensure that every relaxed Dirichlet problem
on $\Om$ can be
approximated in a convenient sense by classical Dirichlet problems on a
suitable sequence of open
sets $(\Ome)$ included in $\Om$.
}

\goodbreak\bigskip
\noindent
{\it Reconstructing the measure $\mu$\/}\nobreak

We now want to reconstruct the measure $\mu$ from one particular
solution of the relaxed Dirichlet problem (\rf{bc1}). In view of the
applications we consider also solutions of the equation in (\rf{bc1})
which do not necessarily satisfy the homogeneous
Dirichlet boundary condition on $\partial\Om$, but we study
only the case where the solution and the right hand side are
nonnegative. Let us fix
$$
A\in\mat\,,\qquad
\mu\in\mo\,,\qquad
\lambda\in\hmp\,,
\leqno(\frm{be0})
$$
and a solution $\om$ to the problem
$$
\cases{\om\in\X\,,&
\cr\cr
\displaystyle
\into A\, D\om\, Dy\,dx +
\into \om\, y \,d\mu
= \into  y \,d\lambda \qquad \forall y\in\x\,,&
\cr}\leqno(\frm{be1})
$$
which satisfies
$$
\om\ge 0\quad\hbox{q.e.\ in }\,\Om \,.
\leqno(\frm{bg10})
$$
\rem{\thm{bi02}}{From the Lax-Milgram lemma, there exists a solution
of (\rf{be1}) which belongs to $\ho$; by the
comparison principle (Theorem~2.10 in
[\rf{DM-Mos-86}]) this solution satisfies (\rf{bg10}), so that the set
of such functions $\om$ is not empty.
}
\bigskip

The following proposition (proved in [\rf{DM-Mos-86}],
Proposition~2.6) will be frequently used throughout the paper.

\prop{\thm{bc03}}{Assume (\rf{be0}), (\rf{be1}), and (\rf{bg10}). Then
there exists
$\nu\in\hmp$ such that
$$
-\div(A\, D\om) + \nu=\lambda \quad \hbox{in }\, \distr\,.
\leqno(\frm{bc5})
$$
}

For technical reasons, the reconstruction of the measure $\mu$ from
$\om$ requires the following assumption: for every quasi open set
$U$ in $\Om$ we have
$$
\cp (U\cap\{\om=0\})>0\quad \Longrightarrow \quad\lambda(U)>0 \,.
\leqno(\frm{be2})
$$
\negbigskip

\rem{\thm{bg02}}{Condition (\rf{be2}) is satisfied in
the following (extreme) cases:
 \smallskip
 \item{(a)}$\om>0$ q.e.\ in $\Om$;
 \smallskip
 \item{(b)} $\lambda(U)>0$ for every quasi open set $U\subseteq\Om$
 with $\cp(U)>0$.
 \smallskip
 \noindent Note that (b) is always satisfied if $\lambda(U)=\intu
 f\,dx$ with $f\in \loneloc$
 and $f>0$ a.e.\ in $\Om$, since, by Lemma \rf{be1.1} and (\rf{bg60}),
 every quasi open
 set with positive capacity has positive Lebesgue measure.
}

\prop{\thm{be01}}{Assume (\rf{be0}), (\rf{be1}), (\rf{bg10}), and
(\rf{be2}).
Then
$$
u\in\x \quad \Longrightarrow \quad  u=0
\quad \hbox{q.e.\ in }\,\{\om=0\}\,.
\leqno(\frm{bg20})
$$
Moreover for every Borel set $B\subseteq \Om$
$$
{\cp(B\cap\{{\om=0}\})>0} \quad \Longrightarrow \quad \mu(B)={+}\infty
\,.\leqno(\frm{bg21})
$$
}
\negbigskip

\proof The proof is along the lines of Lemma 3.2 of [\rf{DM-Gar-94}],
with some important variants, due to the fact that now $\lambda$
is not the Lebesgue measure.

To prove (\rf{bg20}) it is enough to consider
a function $u\in\x$ such that $0\le u\le1$ q.e.\ in $\Om$. For every
$k\in\n$
let $u_k$ be the solution of the relaxed Dirichlet problem
$$
\cases{u_{k}\in\x\,,&
\cr\cr
\displaystyle
\into A\, Du_{k} Dy\,dx +
\into u_{k} y \,d\mu +
k \into u_k y\,d\lambda
= k \into  u\, y \,d\lambda
\cr\cr
 \forall y\in\x\,.&
\cr}\leqno(\frm{be3})
$$
By the comparison principle (see [\rf{DM-Mos-86}], Proposition~2.10)
we have $0\le u_k\le k\,\om$ q.e.\ in $\Om$, hence
$u_k=0$ q.e.\ in $\{{\om=0}\}$.

Taking $y=u_{k}-u$ as test function in (\rf{be3}), from (\rf{za12})
we obtain, by using Cauchy inequality,
$$
\displayliness{
\alpha \into |Du_{k}|^{2} dx + \into |u_{k}|^{2} d\mu + 2k \into
|u_{k}-u|^{2} d\lambda \le
\cr
\le  {1\over\alpha}
\into |A\,Du|^{2} dx + \into |u|^{2} d\mu \,.
\cr}
$$
It follows that $(u_{k})$ is bounded in $\ho$ and converges to $u$
strongly in $L^{2}(\Om,\lambda)$.
Therefore a
subsequence, still denoted by $(u_{k})$, converges weakly in $\ho$ to
some function $v$ in $\ho$ such that $v=u$ $\lambda$-a.e.\ in $\Om$.
Since $u_k=0$ q.e.\ in $\{{\om=0}\}$, and since suitable convex
combinations of $(u_k)$ converge to $v$ strongly in $\ho$,
we conclude that $v=0$ q.e.\ in $\{{\om=0}\}$. Let $V=\{{v\neq u}\}$.
Then $V$ is quasi open and $\lambda(V)=0$. It follows from (\rf{be2})
that
$\cp(V\cap\{{\om=0}\})=0$. As $u=v$ in ${\Om\setminus V}$ and
$v=0$ q.e.\ in $\{{\om=0}\}$,
this implies that $u=0$ q.e.\ in $\{{\om=0}\}$.

Let is prove (\rf{bg21}). Let $U$ be a quasi open subset of $\Om$
such that $\mu(U)<{+}\infty$.
By Lemma \rf{be1.1} there exists an increasing sequence $(z_k)$ in $\ho$
converging to $1_U$ pointwise q.e.\ in $\Om$ and such that $0\le z_k\le
1_U$
q.e.\ in $\Om$ for every $k\in\n$. As $\mu(U)<{+}\infty$,
each function $z_k$
belongs to $\lmu$, hence $z_k=0$ q.e.\ on $\{{\om=0}\}$ by the
previous step. This implies that $1_U=0$ q.e.\ on $\{{\om=0}\}$, hence
${\cp(U\cap\{{\om=0}\})=0}$.

Let us consider a Borel set $B$ with ${\cp(B\cap\{{\om=0}\})>0}$.
For
every quasi open set $U$ containing $B$ we have
${\cp(U\cap\{{\om=0}\})>0}$, hence $\mu(U)={+}\infty$ by the
previous
step of the proof. Then the regularity property (b) in the definition of
$\mo$
implies that $\mu(B)={+}\infty$.
\endproof

\prop{\thm{bf01}}{Assume (\rf{be0}), (\rf{be1}), (\rf{bg10}), and
(\rf{be2}),
 and let $\nu$ be
the measure of $\hmp$ defined in (\rf{bc5}). Then for
every Borel set $B\subseteq\Om$ we have
$$
\mu(B)=\cases{\displaystyle\int_B
{d\nu\over \om}\,,& if $\cp(B\cap\{\om=0\})=0$,\cr
\cr
+\infty\,,&if $\cp(B\cap\{\om=0\})>0$,\cr}
\leqno(\frm{bg13})
$$
and
$$
\nu(B\cap \{\om>0\})=\int_{B}\om\,d\mu\,.
\leqno(\frm{bg14})
$$
In particular, this implies that $\nu=\om\mu$ on ${\{\om>0\}}$.
}

\proof The proof is along the lines of Lemma 3.3 and Proposition 3.4 of
[\rf{DM-Gar-94}].
For every $\eta>0$ let $\nu_{\eta}$ be the Borel measure
defined by
$$
\nu_{\eta}(B)=\int_{B\cap \{\om>\eta\}}\om\,d\mu\,.
\leqno(\frm{bg29})
$$
As $\om\in\lmu$, we have
$\nu_{\eta}(\Om)\le {1\over\eta} \into \om^{2} d\mu <+\infty$.
Let us prove that
$$
\nu_{\eta}(B)=\nu(B\cap \{\om>\eta\})\,,
\leqno(\frm{bg30})
$$
for every Borel set $B\subseteq\Om$. Since $\nu_{\eta}$ is a Radon
measure, it is enough to prove that $\nu_{\eta}(U)=\nu(U\cap
\{{\om>\eta}\})$ for every open set $U\subseteq\Om$. Let us fix an
open set $U$, and
let $U_{\eta}=U\cap \{{\om>\eta}\}$. As $U_{\eta}$ is quasi open, by
Lemma \rf{be1.1} there exists an increasing sequence $(z_{k})$ of
nonnegative functions of $\ho$ converging to $1^{}_{U_{\eta}}$
pointwise
q.e.\ in~$\Om$. Since $\mu(U_{\eta})<+\infty$, the functions $z_{k}$
belong to $\lmu$. Using $z_{k}$ as test function in (\rf{be1}) and
(\rf{bc5}) we obtain
$$
\into z_{k} d\nu = \into \om\, z_{k} d\mu \,.
$$
Taking the limit as $k$ tends to $\infty$ we get
$\nu(U\cap \{{\om>\eta}\})=\nu_{\eta}(U_{\eta})=\nu_{\eta}(U)$,
which proves (\rf{bg30}).
When $\eta$ tends to $0$, we obtain (\rf{bg14}) from
(\rf{bg29}) and (\rf{bg30}) (recall that $\om\ge0$ q.e.\ in $\Om$).

{}From (\rf{bg14}) we have
$$
\mu(B\cap\{\om>\eta\}) = \int_{B\cap\{\om>\eta\}} {d\nu\over \om} \,,
$$
for every Borel set $B\subseteq\Om$ and every $\eta>0$. Taking the limit
as
$\eta$ tends to $0$ we obtain
$$
\mu(B)= \int_{B} {d\nu\over \om}\,,
\leqno(\frm{ma2})
$$
for every Borel set $B\subseteq{\{\om>0\}}$. Since $\mu$ vanishes on
all sets with capacity zero, (\rf{ma2}) holds also when
$\cp(B\cap\{\om=0\})=0$. Finally,
if $\cp(B\cap\{\om=0\})>0$, then $\mu(B)=+\infty$ by
Proposition~\rf{be01}.
\endproof

\goodbreak\bigskip
\noindent{\it Density and uniqueness results\/}\nobreak

In the next proposition we assume, in addition, that
$$
\om\in\linfty\,.
\leqno(\frm{ma1})
$$
The following density result will be crucial in Sections
\rf{corrector} and \rf{proofs}. The proof is along the lines of
Proposition 5.5 of [\rf{DM-Mur-97}], with one important variant, due
to the fact that now the solutions $u_{k}$ of the penalized problem
(\rf{be3}) may not converge to $u$ weakly in $\ho$ (see the proof of
Proposition~\rf{be01}).

\prop{\thm{bf02}}{Assume (\rf{be0}), (\rf{be1}), (\rf{bg10}),
(\rf{be2}),  and  (\rf{ma1}).
Then the set
$\{\om\,\varphi:{\varphi\in\cinftyo}\}$ is dense in $\x$.
}

\proof
For every $u\in\x$ we have to construct a
sequence $(\varphi_{k})$ in $\cinftyo$ such that
$(\om\,\varphi_{k})$ converges to $u$ both in $\ho$ and in $\lmu$.
Clearly it is enough to consider the case $u\ge0$ q.e.\ in $\Om$. For
every
$j\in\n$ let $v_j=u\land (j\,\om)$. Since $\om\ge0$ q.e.\ in $\Om$ and
 $u=0$ q.e.\ in $\{\om=0\}$ by Proposition~\rf{be01}, the sequence
 $(v_j)$ is nondecreasing and converges to $u$ q.e.\ in $\Om$. By
 Lemma~1.6 of [\rf{DM-83}] there exists a sequence $(u_j)$ in $\ho$,
 converging to $u$ strongly in $\ho$, such that $0\le u_j\le
 v_j\le u$ q.e.\ in $\Om$ for every $j\in\n$. By the
dominated convergence theorem it turns out that
 $(u_j)$ converges to $u$ in
 $\lmu$ too.

We are thus reconduced to the case where $u\in\ho$ is such that
$0\le u\le c\,\om$ q.e.\ in $\Om$ for some constant $c>0$.
Since $\{({u-c\,\e})^+>0\}\subseteq\{{\om>\e}\}$,
 and $({u-c\,\e})^+$
converges to~$u$ in $\x$ as $\e$ tends to $0$, we may
also assume that there exists $\e>0$ such that
$\{{u>0}\}\subseteq\{{\om>\e}\}$. Then
${u/ \om} = {u/(\om\lor\e)}$. Since $\om\in\h\cap\linfty$, we have
$u\in\ho\cap\linfty$,
and thus  $u/\om\in\ho\cap\linfty$. Therefore there exists
a sequence $(\varphi_{k})$ in $\cinftyo$, bounded in $\linfty$,
which converges to $z=u/\om$ strongly in $\ho$ and q.e.\ in $\Om$,
hence $\mu$-a.e.\ in $\Om$. Since $\om\in\h\cap\linfty$, the
sequence $(\om\,\varphi_{k})$ converges to $\om\,z=u$ strongly
in $\ho$. As $\om\in\lmu$ and $(\varphi_{k})$ is
bounded in $\linftymu$ and converges to $z=u/\om$ $\mu$-a.e.\ in
$\Om$, by the
dominated convergence theorem the sequence $(\om\, \varphi_{k})$
converges to $\om\, z=u$ strongly in $\lmu$.
\endproof

\bigskip
The following uniqueness result will be crucial in
Theorems \rf{bc020} and~\rf{bg04}. The proof is along the lines of
Lemma~3.5 of [\rf{DM-Gar-94}], with one important variant, due to the
fact that now the condition $\into u^{2} d\lambda=0$ does not imply
that $u=0$ q.e.\ in~$\Om$.

\prop{\thm{bf03}}{Assume (\rf{be0}), (\rf{be1}), (\rf{bg10}),
(\rf{be2}), and (\rf{ma1}).
 Let
$u$ be a solution of the problem
$$
\cases{u \in\ho\cap\linfty\,,&
\cr\cr
\displaystyle
\into A \,D\varphi\,D u \, \om \,dx
- \into A \,D\om\, D\varphi \, u \,dx
+ \into u \,\varphi \,d\lambda=0
\qquad \forall \varphi\in \cinftyo \,.&
\cr}\leqno(\frm{bf10})
$$
Then $u=0$ q.e.\ in $\Om$.
}

\proof Since $\om\in\ho\cap\linfty$, it is easy to see that the
equation in (\rf{bf10}) is satisfied also for
$\varphi\in\ho\cap\linfty$. Using $\varphi=u$ as test function in
this equation we obtain
$$
\into A \,Du \,Du\, \om \,dx
- {1\over 2} \into A \,D\om\, D(u^{2}) \,dx
+ \into u^{2} d\lambda=0\,.
\leqno(\frm{bg1})
$$
Using $y=u^{2}$ as test function in (\rf{be1}), from (\rf{bg1})
we get
$$
\into A \,Du \,Du\, \om \,dx
+ {1\over 2} \into \om\, u^{2} d\mu
+ {1\over 2} \into u^{2} d\lambda=0\,.
$$
This implies
$$
\ldisplaylinesno{
Du=0 \quad \hbox{a.e.\ in }\,\{\om>0\} \,,
&(\frm{bg2})
\cr
u=0 \quad \lambda\hbox{-a.e.\ in }\, \Om \,.
&(\frm{bg3})
\cr}
$$

Let $U=\{u\neq 0\}$. Then $U$ is quasi open and $\lambda(U)=0$ by
(\rf{bg3}). Therefore (\rf{be2}) implies that $u=0$ q.e.\ in
$\{{\om=0}\}$, and consequently $Du=0$ a.e.\ in $\{{\om=0}\}$. By
(\rf{bg2}) we conclude that $Du=0$ a.e.\ in $\Om$. Since
$u\in\ho$, this yields $u=0$ q.e.\ in $\Om$.
\endproof

\parag{\chp{ho}}{A global convergence result}

For every $\e\ge0$ we consider a
matrix
$\Ae$ in $\mat$ and a measure
$\mue$  in $\mo$, that will remain fixed
throughout the rest of the paper.
We assume that
$$
\Ae\ \ H\hbox{-converges to }\ \Ao\,.
\leqno(\frm{ab1.5})
$$

In this section we use a duality argument
to prove that, under suitable hypotheses on $(\mue)$
(which are always satisfied by
a subsequence), the solutions $\ue$ of
the relaxed Dirichlet problems (\rf{bc1}) for
$A=\Ae$ and $\mu=\mue$ converge 
to the solution $\uo$ of the relaxed Dirichlet problem
for $A=\Ao$ and $\mu=\muo$.

\goodbreak\bigskip
\noindent
{\it Definition of special test fuctions\/}\nobreak

For every $\e\ge 0$ we define the functions $\we$ and $\wbe$ as
the unique solutions to
the problems
$$
\ldisplaylinesno{
\cases{\we\in\xe\,,&
\cr\cr
\displaystyle
\int_{\Om}\Ae D\we Dy\,dx +
\int_{\Om} \we y \,d\mue
=\int_{\Om}y\,dx \qquad \forall y\in\xe\,,&
\cr}
&(\frm{ab7})
\cr\cr
\cases{\wbe\in\xe\,,&
\cr\cr
\displaystyle
\int_{\Om}\Abe D\wbe Dy\,dx +
\int_{\Om} \wbe y \,d\mue
=\int_{\Om}y\,dx \qquad \forall y\in\xe\,.&
\cr}
&(\frm{ab7'})
}
$$

By the comparison principle (Theorem 2.10 of [\rf{DM-Mos-86}]) we have
$$
\we\ge 0 \quad\hbox{and}\quad \wbe\ge 0
\quad \hbox{q.e.\ in }\,\Om\,.
\leqno(\frm{bh1})
$$
Moreover, by the maximum principle, we have also
$$
\sup_{\e\ge 0} \|\we\|_{\linfty} <+\infty
\qquad\hbox{and}\qquad
\sup_{\e\ge 0} \|\wbe\|_{\linfty} <+\infty
\leqno(\frm{ab8.5})
$$
(see [\rf{DM-Gar-94}], Section~3).
By Proposition \rf{bc03} there exists two measures
$\nue$ and $\nube$ in $\hmp$ such that
$$
-\div(\Ae D\we) + \nue=1\,, \quad
-\div(\Abe D\wbe) + \nube=1
\quad \hbox{in }\, \distr\,.
\leqno(\frm{ab8})
$$
Finally, from (\rf{bc10})
we obtain
$$
\ldisplaylinesno{
\sup_{\e\ge 0} \into |D\we|^{2} dx <+\infty \,, \qquad
\sup_{\e\ge 0} \into |D\wbe|^{2} dx <+\infty \,,
&(\frm{ab8.6})
\cr
\sup_{\e\ge 0} \into |\we|^{2} d\mue <+\infty \,,
\qquad
\sup_{\e\ge 0} \into |\wbe|^{2} d\mue <+\infty \,.
&(\frm{ab8.7})
\cr}
$$

\goodbreak
\bigskip
\noindent
{\it The main convergence result\/}\nobreak

Given, for every $\e\ge 0$,  $\fe$ and $\fbe$ in $\hm$, we consider the
solutions $\ue$ and $\ube$ to the following problems
$$
\ldisplaylinesno{
\cases{\ue\in\xe\,,&
\cr\cr
\displaystyle
\into\Ae D\ue Dy\,dx +
\into \ue y \,d\mue
=\langle \fe, y \rangle \qquad \forall y\in\xe\,,&
\cr}&(\frm{bc30})
\cr\cr
\cases{\ube\in\xe\,,&
\cr\cr
\displaystyle
\into\Abe D\ube Dy\,dx +
\into \ube y \,d\mue
=\langle \fbe, y \rangle \qquad \forall y\in\xe\,.&
\cr}&(\frm{bc31})
\cr}
$$

\th{\thm{bc020}}{Assume (\rf{ab1.5}) and let $\we$ and $\wbe$ be the
solutions of (\rf{ab7}) and (\rf{ab7'}). The following conditions are
equivalent:
\smallskip
\item{\rm (a)}$\we\wto \wo$ weakly in $\ho$;
\smallskip
\item{\rm (b)}$\wbe\wto \wbo$ weakly in $\ho$;
\smallskip
\item{\rm (c)} for every $(\fe)$ and $(\ue)$ satisfying (\rf{bc30}),
if $\fe\to \fo$ strongly in $\hm$, then
$\ue\wto \uo$ weakly in $\ho$;
\smallskip
\item{\rm (d)} for every $(\fbe)$ and $(\ube)$ satisfying (\rf{bc31}),
if $\fbe\to \fbo$ strongly in $\hm$, then
$\ube\wto \ubo$ weakly in $\ho$.
\par
}

\proof ${\rm (a)}\Rightarrow {\rm (d)}$. Assume (a). By (\rf{bc10})
it is enough
to prove (d) when $\fbe=\fbo=\fb\in\linfty$. Since the equation is
linear,
it suffices to consider the case $0\le \fb\le 1$ a.e.\ in $\Om$, so that

$0\le\ube\le\wbe$ q.e.\ in $\Om$ by the comparison principle
(Theorem 2.10 of [\rf{DM-Mos-86}]).

By (\rf{bc10}) the sequence $(\ube)$ is bounded in $\ho$ and by
(\rf{ab8.5}) it is bounded in $\linfty$. Extracting a
subsequence, we may assume that
$$
\ube \wto \ubs \qquad \hbox{weakly in }\, \ho \,,
\leqno(\frm{bc40})
$$
for some function $\ubs\in\ho\cap\linfty$. We want to show that
$\ubs=\ubo$.
Since the limit does not depend on the
subsequence, this will prove that the whole sequence $(\ube)$ converges
to $\ubo$.

By Proposition \rf{bc03} we have
$$
-\div(\Abe D\ube) + \gabe=\fb
\quad \hbox{in }\, \distr\,,
\leqno(\frm{ma3})
$$
for some $\gabe\in\hmp$.
By Theorem \rf{za006}, from (\rf{ab8}) and (\rf{ma3}) we deduce that
$$
\displayliness{
\Ae D\we \wto \Ao D\wo \qquad\hbox{weakly in }\,\ltwon\,,
\cr
\Abe D\ube \wto \Abo D\ubs \qquad\hbox{weakly in }\,\ltwon\,.
\cr}\leqno(\frm{ab11.2})
$$
Let $\varphi\in\cinftyo$.
Using $y=\we\varphi$ as test function in (\rf{bc31}) and
$y=\ube\varphi$ as test function in (\rf{ab7}), by difference we
obtain
$$
\displayliness{
\into \Abe D\ube D\varphi\, \we dx
- \into \Ae D\we D\varphi \,\ube dx =
\cr
= \into \fb\,\we \varphi \,dx - \into \ube \varphi \,dx \,,
\cr}\leqno(\frm{bc41})
$$
for every $\e\ge 0$.
Since $(\we)$ converges to $\wo$ strongly in $\ltwo$ by (a) and
$(\ube)$ converges to $\ubs$ strongly in $\ltwo$ by (\rf{bc40}),
using (\rf{ab11.2}) we can pass to the
limit in each term of (\rf{bc41}) and we obtain
$$
\displayliness{
\into \Abo D\ubs \,D\varphi \,\wo dx
- \into \Ao D\wo D\varphi \,\ubs \,dx =
\cr
= \into \fb \,\wo \varphi \,dx - \into \ubs \,\varphi \,dx \,.
\cr}\leqno(\frm{bc42})
$$
Since (\rf{bc41}), with $\e=0$, and (\rf{bc42})
hold for every $\varphi\in\cinftyo$, the difference $u=\ubo-\ubs$
belongs to $\ho\cap\linfty$ and
satisfies (\rf{bf10}) with $A=\Ao$, $\om=\wo$, and $\lambda=1$.
This implies $\ubs=\ubo$ q.e.\ in $\Om$
by Proposition \rf{bf03}.

\bigskip
\noindent ${\rm (d)}\Rightarrow {\rm (b)}$.  It is enough to take
$\fbe=\fbo=1$ in condition~(d).

\bigskip
\noindent
${\rm (b)}\Rightarrow {\rm (c)}$. Since $(\Abe)$ $H$-converges to
$\Abo$,
we can replace $\Ae$ by $\Abe$ and $\fbe$ by $\fe$
in the proof of  the implication $(a)\Rightarrow (d)$.

\bigskip
\noindent ${\rm (c)}\Rightarrow {\rm (a)}$.
It is enough to take $\fe=\fo=1$ in
condition~(c).
\endproof

\goodbreak\bigskip
\noindent
{\it A compactness result\/}\nobreak

We now prove that the conditions of Theorem~\rf{bc020} are always
satisfied by a subsequence.

\th{\thm{bd01}}{Assume (\rf{ab1.5}).
For every sequence $(\mue)_{\e>0}$ in $\mo$ there
exist a subsequence, still denoted by $(\mue)$, and a measure
$\muo$ in $\mo$, such that
the equivalent conditions (a)--(d) of Theorem~\rf{bc020} are satisfied.
}

\proof By (\rf{ab8.6}) the sequence $(\we)$ is bounded in $\ho$.
Passing to a subsequence, we may assume that
$(\we)$ converges weakly in $\ho$ to
some function $\ws\in\ho$. By (\rf{bh1}) we have $\ws\ge 0$ q.e.\ in
$\Om$. Now we want to construct a measure
$\muo\in\mo$ such that $\ws$ coincides with the solution
$\wo$ of (\rf{ab7}) for $\e=0$.

By (\rf{ab8}) and Theorem \rf{za006} the sequence
$(\Ae D\we)$ converges to $\Ao D\ws$ weakly in $\ltwon$.
Therefore $(\nue)$ converges to $\nus$ weakly in $\hm$,
where $\nus\in \hmp$ is defined by
$$
-\div(\Ao D\ws) + \nus=1\quad \hbox{in }\, \distr\,.
\leqno(\frm{ab12})
$$

Let us define the measure $\muo$ by
$$
\muo(B)=
\cases{\displaystyle \int_{B} {d\nus\over\ws} & if
$\cp(B\cap\{\ws=0\})=0$,
\cr\cr
+\infty& if
$\cp(B\cap\{\ws=0\})>0$.
\cr}\leqno(\frm{ab13})
$$
Using (\rf{ab12}), from Proposition~3.4 of [\rf{DM-Gar-94}] we obtain
that
$\muo\in\mo$ and that $\ws$ coincides with the unique solution $\wo$ to
problem (\rf{ab7}) for $\e=0$. This shows that condition~(a) of
Theorem~\rf{bc020} is satisfied.
\endproof


\goodbreak\bigskip
\noindent{\it More general test functions\/}\nobreak

We introduce now a more general family of test functions $(\ome)$.
While it is very difficult to compute explicitly the functions $\we$
defined by (\rf{ab7}), in some interesting situations it will be very
easy
to construct explicitly the new family $(\ome)$, from which one can
determine immediately the limit measure $\muo$.

For every $\e\ge 0$ let $\lae\in\hmp$ and
let $\ome$ be a solution of the problem
$$
\cases{
\ome\in \Xe\,,
&
\cr\cr
\displaystyle \int_{\Om}\Ae D\ome Dy\,dx +
\int_{\Om} \ome y \,d\mue
=\int_{\Om} y\,d\lae \qquad \forall y\in\xe\,.
&
\cr}
\leqno(\frm{ba1})
$$
We assume that
$$
\ldisplaylinesno{
\lae\in\hmp \quad \hbox{for every }\, \e\ge 0 \,,
&(\frm{ba2})
\cr
\lae\to\lao\qquad\hbox{strongly in }\,\hm \,,
&(\frm{ba2.5})
\cr
\ome\ge 0 \quad\hbox{q.e.\ in }\,\Om \,\hbox{ for every }\, \e\ge 0 \,,
&(\frm{ba3})
\cr
\ome\wto\omo\qquad\hbox{weakly in }\,\h \,.
&(\frm{ba4})
\cr}
$$
Moreover we assume that for every quasi open set
$U$ in $\Om$ we have
$$
\cp (U\cap\{\omo=0\})>0\quad \Longrightarrow \quad\lao(U)>0 \,,
\leqno(\frm{bg51})
$$
and that
$$
\omo\in\linfty\,.
\leqno(\frm{bg50})
$$

\rem{\thm{bg03}}{If condition (a) of Theorem~\rf{bc020} is satisfied,
then the functions $\we$, $\e\ge0$, defined by (\rf{ab7})
satisfy conditions (\rf{ba1})--(\rf{bg50}) with $\lae=1$ for every
$\e\ge0$
(see Remark~\rf{bg02}). Other sequences $(\ome)_{\e\ge0}$ and
$(\lae)_{\e\ge0}$ satisfying (\rf{ba1})--(\rf{bg50}), with
$\omo=1$, are constructed in [\rf{Cio-Mur}] 
when $\muo\in\hmp$.

If conditions (\rf{ba1})--(\rf{bg50}) are satisfied in $\Om$, then
they are satisfied in every open set $U\subseteq\Om$.
}

\th{\thm{bg04}}{Assume that (\rf{ab1.5}) holds and that
$(\ome)_{\e\ge0}$
and $(\lae)_{\e\ge0}$ satisfy (\rf{ba1})--(\rf{bg50}). Then the
equivalent
conditions (a)--(d) of Theorem \rf{bc020} are fulfilled.
}

\proof We will prove that condition (b) holds.
By (\rf{ab8.6}) the sequence $(\wbe)$ is bounded
in $\ho$ and by (\rf{ab8.5}) it is bounded in $\linfty$. Extracting a
subsequence, we may assume that
$$
\wbe \wto \wbs \qquad \hbox{weakly in }\, \ho \,,
\leqno(\frm{bc400})
$$
for some function $\wbs\in\ho\cap\linfty$. We will show that
$\wbs=\wbo$.
Since the limit does not depend on the
subsequence, this will prove that the whole sequence $(\wbe)$ converges
to $\wbo$.

By Proposition \rf{bc03} and
Theorem \rf{za006} we have
$$
\displayliness{
\Ae D\ome \wto \Ao D\omo \qquad\hbox{weakly in }\,\ltwon\,,
\cr
\Abe D\wbe \wto \Abo D\wbo \qquad\hbox{weakly in }\,\ltwon\,.
\cr}\leqno(\frm{ab11.20})
$$
Let $\varphi\in\cinftyo$.
Using $y=\ome\varphi$ as test function in (\rf{ab7'}) and
$y=\wbe\varphi$ as test function in (\rf{ba1}),
by difference we
obtain
$$
\displayliness{
\into \Abe D\wbe D\varphi \,\ome dx
- \into \Ae D\ome D\varphi \,\wbe dx =
\cr
= \into \ome \varphi \,dx - \into \wbe \varphi \,d\lae \,,
\cr}\leqno(\frm{bc410})
$$
for every $\e\ge 0$.
Since $(\ome)$ converges to $\omo$ strongly in $\ltwoloc$ by (\rf{ba4})
and $(\wbe)$ converges to $\wbs$  strongly in $\ltwo$ by (\rf{bc400}),
using (\rf{ab11.20}) we can pass to the limit in each term of
(\rf{bc410})
and we obtain
$$
\displayliness{
\into \Abo D\wbs\, D\varphi \,\omo dx
- \into \Ao D\omo D\varphi \,\wbs\, dx =
\cr
= \into \omo \varphi \,dx - \into \wbs\, \varphi \,d\lao \,.
\cr}\leqno(\frm{bc420})
$$
Since (\rf{bc410}), with $\e=0$, and (\rf{bc420})
hold for every $\varphi\in\cinftyo$, the difference $\wbo-\wbs$
belongs to $\ho\cap\linfty$ and
satisfies (\rf{bf10}) with $A=\Ao$ and $\lambda=\lao$.
This implies $\wbs=\wbo$ q.e.\ in $\Om$
by Proposition \rf{bf03}.
\endproof

\goodbreak\bigskip
\noindent
{\it Dirichlet problems on varying domains\/}\nobreak

We conclude this section by considering the particular case of
classical Dirichlet problems on varying domains.
Let $(\Ome)_{\e>0}$ be a sequence of open sets, with
$\Ome\subseteq\Om$, and let $\muo$ be a
measure in $\mo$. For every $\e>0$ let $\we$ and $\wbe$ be
the unique solutions to the problems
$$
\ldisplaylinesno{
\cases{\we\in H^{1}_{0}(\Ome)\,,&
\cr\cr
-\div(\Ae D\we ) = 1 \quad \hbox{in }\, {\cal D}'(\Ome) \,,&
\cr}
&(\frm{ab7a})
\cr\cr
\cases{\wbe\in H^{1}_{0}(\Ome)\,,&
\cr\cr
-\div(\Abe D\wbe ) = 1 \quad \hbox{in }\, {\cal D}'(\Ome) \,,&
\cr}
&(\frm{ab7'a})
}
$$
and let $\wo$ and $\wbo$ be the solutions of (\rf{ab7}) and (\rf{ab7'})
with $\e=0$.

Given $\fe$ and $\fbe$ in $\hm$, for $\e>0$, we consider the
solutions $\ue$ and $\ube$ to the following problems
$$
\ldisplaylinesno{
\cases{\ue\in H^{1}_{0}(\Ome)\,,&
\cr\cr
-\div(\Ae D\ue ) = \fe \quad \hbox{in }\, {\cal D}'(\Ome) \,,&
\cr}&(\frm{bc30a})
\cr\cr
\cases{\ube\in H^{1}_{0}(\Ome)\,,&
\cr\cr
-\div(\Abe D\ube ) = \fbe \quad \hbox{in }\, {\cal D}'(\Ome) \,.&
\cr}&(\frm{bc31a})
\cr}
$$
Given $\fo$ and $\fbo$ in $\hm$, let $\uo$ and $\ubo$ be the solutions
of
(\rf{bc30}) and (\rf{bc31}) with $\e=0$.
All functions in $H^{1}_{0}(\Ome)$ are considered as functions in
$\ho$ which are equal to $0$ q.e.\ in ${\Om\setminus\Ome}$.
(Observe that $\ue$, $\ube$, $\fe$, and $\fbe$ are defined in the
whole of $\Om$, for $\e\ge0$.)

\cor{\thm{bi01}}{Assume (\rf{ab1.5}) and let $\we$ and $\wbe$ be the
solutions of (\rf{ab7a}) and (\rf{ab7'a}) for $\e>0$,
and of (\rf{ab7}) and (\rf{ab7'}) for $\e=0$. The following conditions
are
equivalent:
\smallskip
\item{\rm (a)}$\we\wto \wo$ weakly in $\ho$;
\smallskip
\item{\rm (b)}$\wbe\wto \wbo$ weakly in $\ho$;
\smallskip
\item{\rm (c)} for every $(\fe)$ and $(\ue)$ satisfying (\rf{bc30a})
for $\e>0$ and (\rf{bc30}) for $\e=0$,
if $\fe\to \fo$ strongly in $\hm$, then
$\ue\wto \uo$ weakly in $\ho$;
\smallskip
\item{\rm (d)} for every $(\fbe)$ and $(\ube)$ satisfying (\rf{bc31a})
for $\e>0$ and (\rf{bc31}) for $\e=0$,
if $\fbe\to \fbo$ strongly in $\hm$, then
$\ube\wto \ubo$ weakly in $\ho$.
\par
}

\proof For every $\e>0$ let $\mu_{\Ome}$ be the measures introduced
in (\rf{bb1}) with $U=\Ome$. By Remark \rf{bc01} the functions $\we$
and $\wbe$ defined in (\rf{ab7a}) and (\rf{ab7'a}) coincide with the
solutions of (\rf{ab7}) and (\rf{ab7'}) with $\mue=\mu_{\Ome}$.
For the same reason the functions $\ue$
and $\ube$ defined in (\rf{bc30a}) and (\rf{bc31a}) coincide with the
solutions of (\rf{bc30}) and (\rf{bc31}) with $\mue=\mu_{\Ome}$. The
conclusion follows now from Theorem~\rf{bc020}.
\endproof

\rem{\thm{bi02}}{Let $(\lae)$ be a sequence in $\hmp$ and, for
every $\e>0$, let $\ome$ be a function in $\h$ such that $\ome=0$
q.e.\ in ${\Om\setminus\Ome}$ and
$$
-\div(\Ae D\ome ) = \lae \quad \hbox{in }\, {\cal D}'(\Ome) \,.
$$
Let $\lao\in\hmp$ and let $\omo$ be a solution of (\rf{ba1}) with
$\e=0$. If conditions (\rf{ba2})--(\rf{bg50}) are satisfied, then the
equivalent conditions (a)--(d) of Corollary \rf{bi01} are satisfied.
To prove this fact, it
is enough to use Remark \rf{bc01} and Theorem~\rf{bg04}.
}

\parag{\chp{examples}}{An example}

In this section we apply Corollary \rf{bi01} and Remark \rf{bi02} to
a model problem that has not yet been considered in the literature.
The purpose of this example is to show that the measure $\muo$ which
appears in the limit problem depends not only on the sequence
$(\Ome)$ and on $\Ao$, but also on the sequence $(\Ae)$.

To simplify the exposition, we assume $n\ge 3$ (the case $n=2$
requires obvious modifications, as in [\rf{Cio-Mur}]). Let us fix an
exponent $\gamma$ with
$$
1<\gamma<{n\over n-2}\,.
\leqno(\frm{md1})
$$
For every $\e>0$ and $i\in\z^n$ we consider the point $\xie=\e i$,
the open ball $\Bie$ with centre $\xie$ and radius $\e^\gamma$, and
the concentric closed ball $\Cie$ with radius $\e^{n\over n-2}$. By
(\rf{md1}) we have $\Cie\subseteq \Bie$ for $0<\e<1$, and the sets
$(\Bie)^{}_{i\in\z^n}$ are pairwise disjoint for $0<\e<2^{1\over
1-\gamma}$.
Given a bounded open set $\Om\subseteq\rn$ we define
$$
\Be=\Om\cap\bigcup_{i\in\z^n} \Bie\,,
\qquad
\Ce=\Om\cap\bigcup_{i\in\z^n} \Cie\,.
$$

Let us fix two constants $a,\,b\in[\alpha,\beta]$ and let us define
the matrices $\Ae$, for $\e\ge 0$, by
$$
\Ae(x)=\cases{a\,I&for $x\in\Om\setminus\Be$,
\cr\cr
b\,I&for $x\in \Be$,
\cr}
\leqno(\frm{mc7})
$$
where we set $B^0=\emptyset$, so that $\Ao(x)=a\,I$ for every $x\in\Om$.

Since $(\Ae)$ converges in measure to $\Ao$ by (\rf{md1}), it is easy
to prove that $(\Ae)$ $H$-converges to $\Ao$.

Finally, let $\Ome= {\Om \setminus \Ce}$ for every $\e>0$.

We will determine $\muo\in\mo$ such that the equivalent conditions
(a)--(d) of Corollary \rf{bi01} are satisfied. Using Remark \rf{bi02}
we will construct, for $\e\ge0$, a measure $\lae$ in $\hmp$ and,
for $\e>0$, a function $\ome$ in $\h$ such that $\ome=0$ q.e.\ in
$\Ce$ and
$$
-\div(\Ae D\ome)=\lae\quad\hbox{in }\,{\cal D}'(\Ome) \,.
\leqno(\frm{me1})
$$
Then we will prove that conditions (\rf{ba2})--(\rf{bg50}) are
satisfied, where $\omo$ is a solution of (\rf{ba1}) with $\e=0$.

For every $\e>0$ and $i\in\z^n$ let
$\omie\in H^1({\Bie\setminus \Cie})$
be the solution of the equation $\Delta\omie=0$
on ${\Bie\setminus
\Cie}$ which satisfies the boundary conditions $\omie=0$ on
$\partial \Cie$ and $\omie=1$ on $\partial\Bie$. By explicit
computation we find that
$$
\omie(x)=
\ce- \ce\e^n |x-\xie|^{2-n}\quad \hbox{for }
\, x\in\Bie\setminus\Cie\,,
\leqno(\frm{mg1})
$$
where
$$
\ce={1\over 1-\e^{n-\gamma(n-2)}}\longrightarrow 1
\leqno(\frm{mg100})
$$
by (\rf{md1}).
For $0<\e<2^{1\over 1-\gamma}$ we define $\ome$ as the function which is
equal to $\omie$ on
$({\Bie\setminus \Cie})\cap\Om$, and is extended by $0$ on $\Ce$ and by
$1$ on ${\Om\setminus\Be}$.

By direct computation we find that
$$
\int_{\Bie\setminus \Cie} |D\omie|^2 dx = (n-2)\, S_{n-1}\, \ce \e^n\,,
$$
where $S_{n-1}$ is the $(n-1)$-dimensional measure of the boundary of
the unit ball in $\rn$. This yields
$$
\into |D\ome|^2 dx \le  (n-2)\, S_{n-1}\,\ce N^\e \e^n\,,
\leqno(\frm{md2})
$$
where $N^\e$ is the number of indices $i\in\z^n$ such that the
distance from $\xie$ to $\Om$ is less than $\e$. Since
$$
\lim_{\e\to 0}\, N^\e\e^n={\rm meas}(\overline\Om)<+\infty\,,
\leqno(\frm{mc3})
$$
from (\rf{mg100}) and (\rf{md2}) we deduce that $(\ome)$ is bounded in
$\h$. As $(\ome)$ converges to $\omo=1$ in measure, we conclude
that $(\ome)$ converges to $\omo$ weakly in $\h$, i.e., condition
(\rf{ba4}) is fulfilled.

Let $\sige$ denote the $(n-1)$-dimensional measure on
$\Om\cap\partial \Be$ and let $\lae$ be the measure defined by
$$
\lae= b\,(n-2)\,\ce \e^{n-\gamma(n-1)} \sige\,.
$$
Since, by (\rf{mg1}),
$$
{\partial \omie\over\partial \nu} = (n-2)\, \ce \e^{n-\gamma(n-1)}
\quad\hbox{on }\,\partial\Bie\,,
$$
we obtain that $-b\,\Delta\ome=\lae$ in ${\cal D}'(\Ome)$. As $D\ome=0$
a.e.\ in ${\Ome\setminus\Be}$, we have $\Ae D\ome = b\, D\ome$
a.e.\ in $\Ome$ by (\rf{mc7}), and we conclude that (\rf{me1}) holds.

{}From the properties of $\sige$ and from (\rf{mg100}) it follows that
$\lae\in\hmp$ and that
$$
\lim_{\e\to0} \into\varphi\, d\lae =
b\,(n-2)\,S_{n-1} \into \varphi\,dx\,,
\leqno(\frm{mc1})
$$
for every $\varphi\in\cinftyo$.

We now define $\muo=\lao=b\,(n-2)\,S_{n-1}$. Then condition (\rf{ba2})
is satisfied and $\omo=1$ is a solution to problem (\rf{ba1}) for
$\e=0$.
Therefore it remains to prove that $(\lae)$ converges to $\lao$
strongly in $\hm$.

To this aim for every $\e>0$ and $i\in\z^n$ we consider the functions
$\vie$ defined by
$$
\vie(x)=\cases{
b\,\ce \e^n |x-\xie|^{2-n} & if $x\in\Die\setminus\Bie$,
\cr\cr
b\,\ce \e^{n-\gamma(n-2)} & if $x\in\Bie$,
\cr}
$$
where $\Die$ is the open ball with centre $\xie$ and radius $\e/2$.
By computing the normal derivatives of $\vie$ on both sides of
$\partial\Bie$ we obtain that
$$
-\Delta \vie=\lae\quad \hbox{on }\, \Die\,,
\leqno(\frm{mc2})
$$
for $0<\e<2^{1\over 1-\gamma}$.

Let $\Eie$ be the open ball with centre $\xie$ and radius $\e/4$. We
take a cut-off function $\phiie\in C^\infty_c(\Die)$ such that
$\phiie=1$ on $\Eie$, and $0\le\phiie\le1$, $|D\phiie|\le c/\e$, and
$|\Delta\phiie|\le c/\e^2$ on $\Die$, where $c$ is a suitable constant
independent of $\e$ and $i$.

Finally, we define $\ve\in\h$ by
$$
\ve=\sum_{i\in\z^n} \phiie\vie\,.
$$
By (\rf{mc2}) we have
$$
-\Delta\ve=\lae+g^\e\,,
\leqno(\frm{mc4})
$$
where
$$
g^\e= -2 \sum_{i\in\z^n} D\phiie D\vie - \sum_{i\in \z^n} \Delta
\phiie \vie \,.
$$

{}From the definition of $\vie$ and from the estimates for $D\phiie$
and $\Delta\phiie$ we obtain that the sequence $(g^\e)$ is bounded
in $\linfty$. Therefore, passing to a subsequence, we may assume that
$$
g^\e \wto g\qquad\hbox{weakly in }\, \ltwo\, \hbox{ and strongly in
}\,\hm\,.
\leqno(\frm{mc5})
$$

Moreover we have, for $0<\e<4^{1\over 1-\gamma}$,
$$
\eqalign{
\into |D\ve|^2 dx &{}\le2  \sum_{i\in\z^n} \big\{
\int_{\Die\setminus \Bie} |D\vie|^2 dx + {1\over\e^2}
\int_{\Die\setminus \Eie} |\vie|^2 dx \big\} \le
\cr
&\le M\,N^\e \e^n (\e^{n-\gamma(n-2)} + \e^2)\,,
\cr}
$$
for a suitable constant $M$ independent of $\e$.
Taking (\rf{md1}) and (\rf{mc3}) into account, we conclude that
$(D\ve)$ converges to $0$ strongly in $\ltwon$, hence $(\Delta \ve)$
converges to $0$ strongly in $\hm$. By (\rf{mc4}) and (\rf{mc5}) we
obtain that $(\lae)$ converges to $-g$ strongly in $\hm$, and by
(\rf{mc1}) we have $-g=b\,(n-2)\,S_{n-1}=\lao$. Since the limit does
not depend on the subsequence, we conclude that $(\lae)$ converges
to $\lao$ strongly in $\hm$.

Therefore, by Remark \rf{bi02}, if $(\fe)$converges to $\fo$ strongly
in $\hm$, then the solutions $\ue$ of the classical Dirichlet problems
$$
\cases{\ue\in H^{1}_{0}(\Ome)\,,&
\cr\cr
-\div(\Ae D\ue ) = \fe \quad \hbox{in }\, {\cal D}'(\Ome) \,,&
\cr}
$$
extended by $0$ on $\Om\setminus\Ome$, converge weakly in $\ho$ to the
solution $\uo$ of the problem
$$
\cases{\uo\in \ho\,,&
\cr\cr
-\div(\Ao D\uo ) +\muo\uo= \fo \quad \hbox{in }\, \distr \,,&
\cr}
$$
where $\muo=b\,(n-2)\,S_{n-1}$.

Note that, if we change the constant $b$ in the definition of $\Ae$
(see (\rf{mc7})), the $H$-limit $\Ao$ does not change, but the
measure $\muo$ changes. This shows that $\muo$ depends on the whole
sequence $(\Ae)$, and not only on $\Ao$.

\parag{\chp{corrector}}{Global and local corrector results}

In this section we prove a corrector result for the solutions of
problems (\rf{bc30}) in the special case $\fe=\fo=f$, with
$f\in\linfty$.
In Section \rf{Le} we shall consider the case where $(\fe)$ converges
to $\fo$ strongly in $\hm$, together with the case of more general data.

Assume that $(\ome)_{\e\ge0}$ and $(\lae)_{\e\ge0}$
satisfy (\rf{ba1})--(\rf{bg50}). In order to obtain the corrector
result we assume, in addition, that
$$
\ldisplaylinesno{
\sup_{\e\ge 0} \|\ome\|_{\linfty} <+\infty \,,
&(\frm{ba5})
\cr
\sup_{\e\ge 0} \into |\ome|^{2} d\mue <+\infty \,.
&(\frm{ba6})
\cr}
$$

\rem{\thm{bg07}}{If conditions
(\rf{ba1})--(\rf{bg50}), (\rf{ba5}), and  (\rf{ba6})
are satisfied in $\Om$, then they are satisfied in every open set
$U\subseteq\Om$.

The functions $\we$ introduced in (\rf{ab7}) satisfy
conditions (\rf{ba5}) and (\rf{ba6}), as stated in (\rf{ab8.5}) and
(\rf{ab8.7}).
}

\goodbreak\bigskip
\noindent
{\it Global corrector result\/}\nobreak

For $j = 1, 2, \ldots, n$ let us fix a sequence
$(\zje)$ in $\h$ satisfying (\rf{ab3})--(\rf{ab5}).
Let $\uo$ be the solution of (\rf{bc30}) with $\e=0$ and
$\fo=f\in\linfty$.
Let us fix $\delta>0$ and $\psid\in\hh$ such that
$$
\beta \into |D\uo-D(\psid\omo)|^{2} dx +
\into |\uo-\psid\omo|^{2} d\muo < \delta \,.
\leqno(\frm{ab17})
$$
Such a $\psid$ exists since the set
$\{\omo\varphi:{\varphi\in\cinftyo}\}$ is dense in $\xo$
by Proposition \rf{bf02}.

For every $\e>0$ let $\ved$ be the function defined by
$$
\ved=(\psid + \sum_{j=1}^{n} \Dj{\psid} \zje) \, \ome\,.
\leqno(\frm{ab18})
$$
By (\rf{ab3}), (\rf{ab4}), (\rf{ba4}), and (\rf{ba5}) we have
$$
\ved \wto \psid\omo   \qquad\hbox{weakly in }\, \h \,
\hbox{ and weakly}^{*} \hbox{ in }\,\linfty\,.
\leqno(\frm{ab18.2})
$$
Moreover we have
$$
D\ved =
(\psid + \sum_{j=1}^{n} \Dj{\psid} \zje) D\ome +
\sum_{j=1}^{n} \Dj{\psid} (\ej + D\zje ) \,\ome
+\sum_{j=1}^{n} D\Dj{\psid} \zje \ome \,.
$$
The last sum in the right hand side converges to $0$ strongly in
$\ltwon$ by
(\rf{ab4}) and (\rf{ba5}), while $(\Dj{\psid} \zje D\ome)$ converges to
$0$
strongly in $\ltwon$ by (\rf{ab4}) and (\rf{ba4}).
Therefore
$$
D\ved  = \psid D\ome + \sum_{j=1}^{n} \Dj{\psid} (\ej + D\zje )\, \ome
+ H_{\delta}^{\e}\,,
\leqno(\frm{ab18.5})
$$
where $(H_{\delta}^{\e})$ converges to $0$ strongly in $\ltwon$ as
$\e$ tends to $0$.

Since $\psid\in W^{1,\infty}(\Om)$ and $(\zje)$ is bounded in
$\linfty$, from (\rf{ba6}) we deduce that
$$
\sup_{\e>0} \into |\ved|^{2} d\mue <+\infty \,.
\leqno(\frm{ab26})
$$

\th{\thm{ab03}}{Assume (\rf{ab1.5}), (\rf{ba1})--(\rf{bg50}),
(\rf{ba5}),
and (\rf{ba6}).
Let $\delta>0$ and
let $\psid$ be a function in $\hh$ which satisfies (\rf{ab17}).
Assume that the functions $\ved$ defined by (\rf{ab18}) belong to
$\ho$. Then for every $f\in\linfty$ the solutions $\ue$ of
problems (\rf{bc30}) with $\fe=f$ satisfy the estimate
$$
\limsup_{\e\to0} \big\{ \alpha \into |D\ue-D\ved|^{2} dx
+ \into |\ue-\ved|^{2} d\mue \big\} < \delta \,.
\leqno(\frm{ab19})
$$
}
\negbigskip
\rem{\thm{ab04}}{In the special case $\uo=\psi\,\omo$,
for some $\psi\in\hh$, we can
take $\psid=\psi$ for every $\delta>0$ in (\rf{ab17}), so that
$$
\ved=\ve= (\psi + \sum_{j=1}^{n} \Dj{\psi} \zje)\, \ome\,.
\leqno(\frm{z1})
$$
Therefore  (\rf{ab19}) implies
$$
\lim_{\e\to0} \big\{ \alpha \into |D\ue-D\ve|^{2} dx
+ \into |\ue-\ve|^{2} d\mue \big\}=0\,,
\leqno(\frm{z2})
$$
which is a corrector result.

When the measures $\mue$ are
fixed and equal to $0$ (so that we can choose $\ome=\omo=1$ and
$\psi=\uo$),  formulas (\rf{z1}) and (\rf{z2})
provide the classical corrector result for $H$-converging operators
stated in (\rf{za25})
(see [\rf{Mur-Tar}] and, in the periodic case, [\rf{BLP}] and
[\rf{Sanchez}]).
When the matrices $\Ae$ are fixed and equal to some
matrix $\Ao$ (so that we can choose $\zje=0$),
formulas (\rf{z1}) and (\rf{z2}) with $\ome=\we$ defined by (\rf{ab7})
provide the corrector result of [\rf{DM-Gar-94}] and
[\rf{DM-Mur-97}]; with a different choice of $\ome$, which leads to
$\omo=1$, the same formulas give also the corrector result of
[\rf{Cio-Mur}] in the periodic case. When
both $\Ae$ and $\mue$ depend on $\e$, but $\omo=1$,
so that we have $\psi=\uo$,
the combination of $H$-converging operators and varying domains
results in the multiplication of the corresponding correctors.

In the general case, $\psid$ and $\ved$ depend on $\delta$ and we
obtain from
(\rf{ab19}) that
$$
D\ue=D\ved+R^{\e}_{\delta}
\qquad\hbox{with}\qquad
\limsup_{\e\to0}\|R^{\e}_{\delta}\|_{\ltwon}^{2} <
{\delta\over\alpha}\,,
$$
which is still a corrector result, but in a more technical form.
}

\goodbreak\bigskip
\noindent{\it Local convergence and corrector results\/}\nobreak

We consider now the case where the functions $\ue$
are solutions of the problems
$$
\cases{
\ue\in \Xe\,,
&
\cr\cr
\displaystyle \int_{\Om}\Ae D\ue Dy\,dx +
\int_{\Om} \ue y \,d\mue
=\into f\, y\,dx \qquad \forall y\in\xe\,,
&
\cr}
\leqno(\frm{ab23})
$$
but are not required to satisfy the
boundary condition $\ue=0$ on $\partial\Om$. We still consider the
case of data $f\in \linfty$. More general data will be studied in
Section~\rf{Le}.

The following theorem is a local version of the convergence result
given in Theorem \rf{bg04}. It will be proved in Section~\rf{proofs}.

\th{\thm{ab061}}{Assume (\rf{ab1.5}), (\rf{ba1})--(\rf{bg50}),
(\rf{ba5}),
and (\rf{ba6}).
Let $f\in\linfty$ and, for every $\e>0$, let $\ue$ be a solution of
(\rf{ab23}).
Assume that
$$
\ue\wto\uo \qquad\hbox{weakly in }\, \h\,,
\leqno(\frm{ab20})
$$
for some function $\uo\in\h$, and that
$$
\ldisplaylinesno{
\sup_{\e> 0} \|\ue\|_{\linfty} <+\infty \,,
&(\frm{bd10})
\cr
\sup_{\e> 0} \into |\ue|^{2} d\mue <+\infty \,.
&(\frm{bd11})
\cr}
$$
Then $\uo$ is a solution of (\rf{ab23}) for $\e=0$.
}
\bigskip
The following lemma, which will be proved in Section~\rf{proofs},
shows that (under the other assumptions of Theorem \rf{ab061})
conditions (\rf{bd10}) and (\rf{bd11}) are always
satisfied in every open set $U\subset\subset\Om$, and also in $\Om$
if every $\ue$ belongs to $\ho$.

\lemma{\thm{ab08.50}}{Assume (\rf{ab1.5}), (\rf{ba1})--(\rf{bg50}),
(\rf{ba5}), and (\rf{ba6}). Let $f\in\linfty$ and, for every $\e>0$,
let $\ue$ be a solution of (\rf{ab23}).
Assume that (\rf{ab20}) holds
for some function $\uo\in\h$. Then
we have
$$
\ldisplaylinesno{
\sup_{\e>0} \|\ue\|_{\linftyu} <+\infty \,,
&(\frm{bd1000})
\cr
\sup_{\e>0} \intu |\ue|^{2} d\mue <+\infty \,,
&(\frm{bd1100})
\cr}
$$
for every open set $U\subset\subset \Om$. If, in addition, $\ue\in\ho$
for every $\e>0$, then (\rf{bd1000}) and (\rf{bd1100}) also hold  for
$U=\Om$.
}

\bigskip

In the next corollary $\hc$ denotes the space of all
functions $u\in \h$ with compact support in $\Om$.
The first assertion of the corollary follows immediately from
Theorem \rf{ab061} and Lemma \rf{ab08.50}, while the last assertion
is easily obtained by approximating any nonnegative function
$y\in\xo$ by the sequence $(\varphi_j\land y)$, where
$\varphi_j\in\cinftyo$ converges to $y$ in $\ho$.

\cor{\thm{ab06}}{Under the assumptions of Lemma \rf{ab08.50},
$\uo$ is a solution to the problem
$$
\cases{
\uo\in\h\cap\lmuoloc\,,
&
\cr\cr
\displaystyle \int_{\Om}\Ao D\uo Dy\,dx +
\int_{\Om} \uo y \,d\muo
=\into f\, y\,dx \qquad \forall y\in\hc\cap\lmuo\,.
&
\cr}
\leqno(\frm{ab24})
$$
If, in addition, $\uo\in\lmuo$, then the last line of (\rf{ab24}) holds
for every $y\in\xo$.
}

\bigskip

Let us fix an open set $U\subset\subset \Om$ and a function
$\zeta\in\cinftyo$ such that
$\zeta=1$ in $U$. Given $\uo\in H^1_{\rm loc}(\Om)\cap\lmuoloc$,
by Proposition \rf{bf02} we can
approximate the function $\zeta\, \uo$ in $\xo$ by functions of the
form $\psi\,\omo$ with $\psi\in\cinftyo$.
Therefore for every $\delta>0$ there exists
$\psid\in \hhu$ such that
$$
\beta \intu |D\uo-D(\psid\omo)|^{2} dx +
\intu |\uo-\psid\omo|^{2} d\muo <  \delta \,.
\leqno(\frm{af51})
$$

The following theorem is a local version of the corrector result
given in Theorem~\rf{ab03}.

\th{\thm{ad01}}{Under the hypotheses of Lemma \rf{ab08.50},
let $U$ be an open set with
$U\subset\subset \Om$,
let $\delta>0$, let $\psid$ be a function in $\hhu$ which satisfies
(\rf{af51}), and let $\ved$ be the functions defined in $U$ by
(\rf{ab18}). Then
$$
\limsup_{\e\to0} \big\{ \alpha \intv  |D\ue-D\ved|^{2} dx
+ \intv  |\ue-\ved|^{2} d\mue \big\}< \delta \,,
\leqno(\frm{ad10})
$$
for every open set $V \subset\subset U$.
}

\bigskip

Theorems~\rf{ab03} and~\rf{ad01} can be deduced from
the following theorem, which will be proved in Section~\rf{proofs}.
Indeed, by Theorem \rf{bg04} and Lemma~\rf{ab08.50},
the assumptions of
Theorem \rf{ab03} imply all assumptions of Theorem \rf{ab061}, so that
(\rf{ab19}) follows from
(\rf{za12}), (\rf{za13}), and (\rf{ab25}) with $\varphi=1$.
Similarly, the assumptions of
Theorem \rf{ad01} imply, by Lemma~\rf{ab08.50}, that
all assumptions of Theorem \rf{ab061} are satisfied in every open
set $U\subset\subset\Om$, so that we can apply Theorem \rf{ab07}
with $\Om$ replaced by $U$ and with $\varphi\in C^\infty_c(U)$
such that
$\varphi=1$ in $V$ and $\varphi\ge 0$ in ${U\setminus V}$.

\th{\thm{ab07}}{Under the hypotheses of Theorem~\rf{ab061},
let $\psi$ be a function in $\hh$, and let $\ve$ be defined by
$$
\ve=(\psi + \sum_{j=1}^{n} \Dj{\psi} \zje) \, \ome\,.
\leqno(\frm{af1})
$$
Then for every $\varphi\in\cinftyo$ we have
$$
\eqalign{
\lim_{\e\to0} &\big\{ \into \Ae D(\ue-\ve)D(\ue-\ve)\,\varphi\, dx
+\into |\ue-\ve|^{2}\varphi \,d\mue \big\} =
\cr
&= \into \Ao D(\uo-\psi\,\omo)D(\uo-\psi\,\omo)\,\varphi\, dx
+\into |\uo-\psi\,\omo|^{2}\varphi \,d\muo\,.
\cr}\leqno(\frm{ab25})
$$
If the functions $\ue$ and $\ve$ belong to $\xe$
for every $\e>0$,  then (\rf{ab25}) also holds with
$\varphi=1$.
}

\parag{\chp{comparison}}{A comparison theorem}

In this section we state and prove a comparison result for the limit
measures $\muoneo$ and $\mutwoo$ corresponding to different sequences
of $H$-convergent matrices $\Aonee$ and $\Atwoe$. This result has its
own interest and will be crucial in the proof of the corrector results
stated in the previous section.

For every $\e\ge 0$ let $\Aonee$ and $\Atwoe$ be two matrices
in $\mat$.
We assume that
$$
\Aie\ \ H\hbox{-converges to }\ \Aio \quad \hbox{for }\,i=1,\,2.
\leqno(\frm{bd3})
$$

For every $\e>0$ let $\mue$ be a measure in $\mo$, and let
$\muoneo$ and $\mutwoo$ be two measures in $\mo$. For $i=1,\,2$ and
$\e>0$ let $\wie$ be the solutions of the problems
$$
\cases{\wie\in\xe\,,&
\cr\cr
\displaystyle
\int_{\Om}\Aie D\wie Dy\,dx +
\int_{\Om} \wie y \,d\mue
=\int_{\Om}y\,dx \qquad \forall y\in\xe\,,&
\cr}
\leqno(\frm{ab7+})
$$
and let $\wio$ be the solutions of the problems
$$
\cases{\wio\in\xio\,,&
\cr\cr
\displaystyle
\int_{\Om}\Aio D\wio Dy\,dx +
\int_{\Om} \wio y \,d\muio
=\int_{\Om}y\,dx \qquad \forall y\in\xio\,.&
\cr}
\leqno(\frm{mx1})
$$
We assume that
$$
\wie\wto\wio\qquad\hbox{weakly in }\,\ho\,.
\leqno(\frm{aa4})
$$
Note that, by Theorem \rf{bd01}, these hypotheses are always
satisfied by a subsequence.

In this section we shall prove the following comparison theorem.

\th{\thm{aa04}}{Assume (\rf{bd3}) and (\rf{aa4}).
Then
$$
\ldisplaylinesno{
{\alpha^{2}\over \beta^{2}} \mutwoo \le
\muoneo \le {\beta^{2}\over \alpha^{2}} \mutwoo \quad\hbox{in }\,\Om\,,
&(\frm{aa29})
\cr
\cp(\{\woneo>0\} \, {\scriptstyle\triangle}\, \{\wtwoo>0\})=0 \,.
&(\frm{aa29+})
\cr}
$$
In particular we have $\lmuoneo=\lmutwoo$.
}
\bigskip
In order to prove Theorem~\rf{aa04}, for $\e\ge 0$ and $i=1,\,2$
we consider the measures $\nuie\in\hmp$ defined by
$$
-\div(\Aie D\wie) + \nuie=1\qquad \hbox{in }\, \distr
\leqno(\frm{aa3})
$$
(see Proposition~\rf{bc03}). By Proposition \rf{bf01} we have
$$
\nuio=\wio \muio \quad \hbox{on}\quad \{\wio>0\}\,.
\leqno(\frm{aa9})
$$

By Theorem~\rf{za006} we have
$$
\Aie D\wie \wto \Aio D\wio \qquad\hbox{weakly in }\,\ltwon \,.
\leqno(\frm{bh3})
$$
Therefore
$$
\nuie\wto\nuio\qquad\hbox{weakly in }\,\hm\,.
\leqno(\frm{aa5})
$$
As $\nuie\ge0$, by Theorem~1 of [\rf{Mur-81}] we have
$$
\psi\,\nuie\to\psi\,\nuio\qquad\hbox{strongly in }\,\hmq\,,
\leqno(\frm{aa6})
$$
for every $\psi\in\cinftyo$ and for every $q<2$.

Let $\zie$ be the solution of the problem
$$
\cases{\zie\in\ho\,,&
\cr\cr
-\div(\Aie D\zie)=-\div(\Aio D\wio)\quad\hbox{in }\,\distr\,.
&\cr}\leqno(\frm{aa11})
$$
By the definition of $H$-convergence we have
$$
\ldisplaylinesno{
\zie\wto\wio\qquad\hbox{weakly in }\,\ho\,,
&(\frm{aa12})
\cr
\Aie D\zie \wto \Aio D\wio
\qquad\hbox{weakly in }\,\ltwon\,.
&(\frm{aa13})
\cr}
$$

\lemma{\thm{aa01}}{For every $\varphi\in\cinftyo$
and $i=1,\,2$ we have
$$
\ldisplaylinesno{
\lim_{\e\to0}\big\{ \into\Aonee D(\wonee-\zonee)
D(\wtwoe-\ztwoe)\varphi \,dx + \into \wonee\wtwoe\varphi \,d\mue
\big\} = \into \wtwoo\varphi \,d\nuoneo \,,
&(\frm{aa14})
\cr
\lim_{\e\to0}\big\{ \into\Aie D(\wie-\zie)
D(\wie-\zie)\varphi \,dx + \into |\wie|^2\varphi \,d\mue
\big\} = \into \wio\varphi \,d\nuio \,.
&(\frm{aa14+})
\cr}
$$
}
\negbigskip


\proof Let us first prove (\rf{aa14}). For every $\e>0$ we write
$$
\into\Aonee D(\wonee-\zonee)
D(\wtwoe-\ztwoe)\varphi \,dx + \into \wonee\wtwoe\varphi \,d\mue
= I^{\e}+I\!I^{\e}+I\!I\!I^{\e}\,,
\leqno(\frm{aa15})
$$
where
$$
\eqalign{I^{\e}&=
\into\Aonee D\wonee
D\wtwoe\varphi \,dx + \into \wonee\wtwoe\varphi \,d\mue\,,
\cr
I\!I^{\e}&=-\into\Aonee D\wonee
D\ztwoe\varphi \,dx \,,
\cr
I\!I\!I^{\e}&=-\into\Aonee D\zonee
D(\wtwoe-\ztwoe)\varphi \,dx\,.
\cr}
$$
Using $y=\wtwoe\varphi$ as test function in (\rf{ab7+}) we get
$$
I^{\e}=\into\wtwoe\varphi\,dx - \into\Aonee D\wonee D\varphi \,\wtwoe
\,dx \,.
$$
Since $(\wtwoe)$ converges to $\wtwoo$
strongly in $\ltwo$ by (\rf{aa4})
and since $(\Aonee D\wonee)$ converges
to $\Aoneo D\woneo$ weakly in
$\ltwon$ by (\rf{bh3}), we have
$$
\displayliness{
\lim_{\e\to0} I^{\e} = \into\wtwoo\varphi\,dx
- \into\Aoneo D\woneo D\varphi \wtwoo \,dx =
\cr
= \into\Aoneo D\woneo
D\wtwoo\varphi \,dx + \into \wtwoo\varphi \,d\nuoneo\,,
\cr}\leqno(\frm{aa17})
$$
where in the last equality we used~(\rf{aa3}) for $\e=0$. Note that we
can not use $\wtwoo\varphi$ as test function in (\rf{mx1}) for $i=1$
because we do not know yet that $\wtwoo\varphi\in \lmuoneo$.

{}From~(\rf{aa3}) we obtain
$$
\displayliness{
I\!I^{\e}= -\into\Aonee D\wonee D(\ztwoe\varphi)\,dx
+ \into\Aonee D\wonee D\varphi \,\ztwoe\,dx=
\cr
=\langle\nuonee,\ztwoe\varphi\rangle
-\into \ztwoe\varphi\,dx
+ \into\Aonee D\wonee D\varphi \,\ztwoe\,dx \,.
\cr}\leqno(\frm{aa18})
$$
Since $(\ztwoe)$ converges to $\wtwoo$ strongly in $\ltwo$
by (\rf{aa12}) and
$(\Aonee D\wonee)$ converges to $\Aoneo D\woneo$ weakly in
$\ltwon$ by (\rf{bh3}), we have
$$
\displayliness{
\lim_{\e\to0} \big\{ -\into \ztwoe\varphi\,dx
+ \into\Aonee D\wonee D\varphi \,\ztwoe\,dx \big\}=
\cr
= -\into \wtwoo\varphi\,dx
+ \into\Aoneo D\woneo D\varphi \,\wtwoo\,dx \,.
\cr}\leqno(\frm{aa19})
$$

We will prove in Lemma~\rf{mx2} that
$$
\lim_{\e\to0}\, \langle \nuonee, \ztwoe\varphi \rangle
= \langle \nuoneo, \wtwoo\varphi \rangle \,.
\leqno(\frm{aa23})
$$
{}From (\rf{aa18}), (\rf{aa19}), and (\rf{aa23}) it follows that
$$
\displayliness{
\lim_{\e\to0}
I\!I^{\e} =
\into\wtwoo\varphi\,d\nuoneo
- \into\wtwoo\varphi\,dx
+ \into\Aoneo D\woneo D\varphi \,\wtwoo\,dx=
\cr
=-\into\Aoneo D\woneo D \wtwoo \varphi\,dx \,,
\cr}\leqno(\frm{aa24})
$$
where the last equality is obtained by using $\wtwoo\varphi$ as
test function in~(\rf{aa3}) for $\e=0$.

{}From (\rf{aa11}) it follows that
$$
\displaylines{I\!I\!I^{\e}= -\into \Aonee D\zonee D((\wtwoe-\ztwoe)\,
\varphi)\, dx
+ \into \Aonee D\zonee D \varphi \,(\wtwoe-\ztwoe)\, dx =
\cr
= -\into \Aoneo D\woneo D((\wtwoe-\ztwoe)\,
\varphi)\, dx
+ \into \Aonee D\zonee D \varphi \,(\wtwoe-\ztwoe)\, dx \,.
\cr}
$$
Since $(\wtwoe-\ztwoe)$ converges to $0$ weakly in $\ho$ and strongly
in $\ltwo$ by (\rf{aa4}) and (\rf{aa12}), while $(\Aonee D\zonee)$
converges to $\Aoneo D\zoneo$
weakly in $\ltwon$ by (\rf{aa13}), we have
$$
\lim_{\e\to0} I\!I\!I^{\e}=0 \,.
\leqno(\frm{aa25})
$$

Equality (\rf{aa14}) now follows from (\rf{aa15}), (\rf{aa17}),
(\rf{aa24}), and (\rf{aa25}).

Let us prove now (\rf{aa14+}) for a given $i=1,\,2$. To this aim for
every $\e>0$ we define $\hat A{}^\e_1=\hat A{}^\e_2=\Aie$ and
$\hat\mu{}^\e=\mue$, so that $\hat w{}^\e_1=\hat w{}^\e_2=\wie$,
$\hat\zeta{}^\e_1=\hat\zeta{}^\e_2=\zie$, and $\hat\nu{}^0_1=
\hat\nu{}^0_2=\nuio$. Applying (\rf{aa14}) in this new setting gives
(\rf{aa14+}).
\endproof

\lemma{\thm{mx2}}{For every $\varphi\in\cinftyo$ we have
$$
\lim_{\e\to0}\, \langle \nuonee, \ztwoe\varphi \rangle
= \langle \nuoneo, \wtwoo\varphi \rangle \,.
\leqno(\frm{aa23+})
$$
}
\negbigskip
\proof
Given $\delta>0$, let $\zeta^{0}\in\cinftyo$ be a function such that
$$
\|\zeta^{0}-\wtwoo\|_{\ho} < \delta \,,
\leqno(\frm{aa20})
$$
and let $\zeta^{\e}$ be the solution of the problem
$$
\cases{\zeta^{\e}\in\ho\,,&
\cr\cr
-\div(\Atwoe D\zeta^{\e})=-\div(\Atwoo D\zeta^{0})\quad
\hbox{in }\,\distr\,.
\cr}\leqno(\frm{aa21})
$$
Using $\zeta^{\e}-\ztwoe$ as test function in (\rf{aa11})
and (\rf{aa21}), from (\rf{za12}) and (\rf{za13}) we obtain
$$
\|\zeta^{\e}-\ztwoe\|_{\ho} \le {\beta\over \alpha}
\|\zeta^{0}-\wtwoo\|_{\ho} \,.
\leqno(\frm{aa22})
$$

As $(\nuonee)$ is bounded in $\hm$, form (\rf{aa20}) and (\rf{aa22})
we obtain that there exists a constant $M$, independent of $\delta$,
such that
$$
\displayliness{
|\langle \nuonee,\ztwoe\varphi\rangle
- \langle\nuoneo, \wtwoo\varphi\rangle |
\le
\cr
\le |\langle \nuonee, (\ztwoe-\zeta^{\e}) \varphi\rangle|
+ |\langle \nuonee, \zeta^{\e} \varphi\rangle
- \langle \nuoneo, \zeta^{0} \varphi\rangle|
+ |\langle \nuoneo, (\zeta^{0}-\wtwoo) \varphi\rangle| \le
\cr
\le M\delta + |\langle \nuonee, \zeta^{\e} \varphi\rangle
- \langle \nuoneo, \zeta^{0} \varphi\rangle|\,.
\cr}\leqno(\frm{bh4})
$$

By Meyers' estimate, there exists $p>2$ such that
$(\zeta^{\e}\varphi)$ is bounded in $\hop$. As $(\zeta^{\e})$
converges to $\zeta^{0}$ weakly in $\ho$ by the definition of
$H$-convergence, we conclude that
$(\zeta^{\e}\varphi)$ converges to $\zeta^{0}\varphi$ weakly in $\hop$.
Since by (\rf{aa6}) the sequence $(\psi\nuonee)$ converges to
$\psi\nuoneo$
strongly in $\hmq$ for $1/p+1/q=1$ and for every $\psi\in\cinftyo$,
we obtain that
$$
\lim_{\e\to0} \,\langle \nuonee, \zeta^{\e}\varphi \rangle
= \lim_{\e\to0} \,\langle \psi \,\nuonee, \zeta^{\e}\varphi \rangle
= \langle \psi \,\nuoneo, \zeta^{0}\varphi \rangle
= \langle \nuoneo, \zeta^{0}\varphi \rangle \,,
$$
where $\psi$ is any function in $\cinftyo$ which is equal to $1$ in a
neighbourhood of $\supp(\varphi)$.

Therefore by (\rf{bh4})
$$
\limsup_{\e\to0}
|\langle \nuonee,\ztwoe\varphi\rangle
- \langle\nuoneo, \wtwoo\varphi\rangle |\le M\delta\,.
$$
As $\delta>0$ is arbitrary, we obtain (\rf{aa23+}).
\endproof

\lemma{\thm{aa02}}{For every $\varphi\in\cinftyo$, with $\varphi\ge0$
in $\Om$, and for every $t>0$ we have
$$
\into\wtwoo\varphi\,d\nuoneo \le {\beta\over\alpha}
\big\{
{t\over 2} \into \woneo\varphi\,d\nuoneo
+ 
{1\over 2t} \into \wtwoo\varphi\,d\nutwoo
\big\}
\,.
\leqno(\frm{aa26})
$$
}
\negbigskip

\proof By (\rf{za12}) and (\rf{za13})  we have the estimates
$$
\displaylines{
\into\Aonee D(\wonee - \zonee) D(\wtwoe - \ztwoe) \varphi \,dx
+ \into \wonee \wtwoe \varphi \,d\mue \le
\cr
\le \beta \into| D(\wonee - \zonee)| | D(\wtwoe - \ztwoe) | \varphi \,dx
+ \into \wonee \wtwoe \varphi \,d\mue \le
\cr
\le {t\over 2} \big\{ \beta \into| D(\wonee - \zonee)|^{2} \varphi \,dx
+ \into |\wonee |^{2} \varphi \,d\mue \big\} +{}
\cr
{}+ {1\over 2t}  \big\{ \beta \into | D(\wtwoe - \ztwoe) | ^{2} \varphi
\,dx
+ \into| \wtwoe |^{2} \varphi \,d\mue \big\} \le
\cr
\le {\beta\over\alpha} {t\over 2} \big\{
\into\Aonee D(\wonee - \zonee) D(\wonee - \zonee) \varphi \,dx
+ \into |\wonee |^{2} \varphi \,d\mue  \big\} +\hbox{}
\cr
\hbox{}+ {\beta\over\alpha} {1\over 2t}
\big\{  \into\Atwoe D(\wtwoe - \ztwoe) D(\wtwoe - \ztwoe) \varphi \,dx
+ \into |\wtwoe|^{2} \varphi \,d\mue  \big\} \,.
\cr}
$$
Inequality (\rf{aa26}) is obtained by applying Lemma~\rf{aa01}.
\endproof

\lemma{\thm{aa03}}{The following inequality holds:
$$
\wtwoo\nuoneo \le {\beta^{2}\over \alpha^{2}}
\woneo\nutwoo \quad \hbox{in }\,\Om \,.
\leqno(\frm{aa27})
$$
}
\negbigskip

\proof Let $\nu=\nuoneo+\nutwoo$. {}From Lemma~\rf{aa02} it follows
that for every $t>0$
$$
\wtwoo {d\nuoneo\over d\nu} \le
{\beta\over\alpha}
\big\{
{t\over 2}
\woneo {d\nuoneo\over d\nu}
+ 
{1\over 2t}
\wtwoo {d\nutwoo\over d\nu}
\big\}
\qquad\nu\hbox{-a.e.\ in }\, \Om\,.
\leqno(\frm{aa28})
$$
If we minimize with respect to $t$ we obtain
$$
\wtwoo {d\nuoneo\over d\nu} \le
{\beta^{2}\over \alpha^{2}}
\woneo {d \nutwoo \over d\nu}
\qquad\nu\hbox{-a.e.\ in }\, \Om\,,
$$
which implies (\rf{aa27}).
\endproof

\proofof{Theorem \rf{aa04}}
We prove only the second inequality in (\rf{aa29}) and
$$
\cp(\{{\wtwoo>0}\}\setminus \{{\woneo>0}\})=0\,.
\leqno(\frm{aa29++})
$$
The other inequality and the equality
$\cp(\{{\woneo>0}\} \setminus \{{\wtwoo>0}\})=0$ are proved by
exchanging the roles of $\Aonee$ and $\Atwoe$.

By (\rf{aa9}) we have $\nutwoo=\wtwoo\mutwoo$ on $\{{\wtwoo>0}\}$, so
that (\rf{aa27}) gives
$$
\nuoneo \le {\beta^{2}\over \alpha^{2}}
\woneo\mutwoo \qquad\hbox{on }\, \{{\wtwoo>0}\}\,.
\leqno(\frm{aa31})
$$
If $y\in\xtwoo$, then $y=0$ q.e.\ on $\{{\wtwoo=0}\}$ (see
Proposition~\rf{be01}). {}From (\rf{aa3}) and (\rf{aa31}) it
follows that
$$
\into \Aoneo D\woneo  Dy\,dx +
{\beta^{2}\over \alpha^{2}}
\into \woneo y \,d\mutwoo
\ge \into y\,dx
\leqno(\frm{aa32})
$$
for every $y\in\xtwoo $ with $y\ge0$ q.e.\ in $\Om$.

Let $w$ be the solution of the problem
$$
\cases{w\in\xtwoo\,,&
\cr\cr
\displaystyle
\into \Aoneo Dw Dy\,dx +
{\beta^{2}\over \alpha^{2}}
\into w y \,d\mutwoo
=\into y\,dx \qquad \forall y\in\xtwoo\,.&
\cr}\leqno(\frm{aa30})
$$
As $0\le (w-\woneo)^{+}\le w$ q.e.\ in $\Om$, the function
$y=(w-\woneo)^{+}$ can be taken as test function in (\rf{aa30}) and
(\rf{aa32}). By difference we obtain
$$
\into \Aoneo D(w-\woneo) D(w-\woneo)^{+} dx +
{\beta^{2}\over \alpha^{2}}
\into (w-\woneo) (w-\woneo)^{+} d\mutwoo \le 0\,,
$$
which implies $(w-\woneo)^{+}=0$ a.e.\ in $\Om$, and hence
$w\le\woneo$ q.e.\ in $\Om$ by (\rf{bg60}).
Therefore
$$
\cp(\{w>0\}\cap\{\woneo=0\})=0 \,.
\leqno(\frm{aa33})
$$

Let us prove that
$$
\cp(\{\wtwoo>0\}\cap\{w=0\})=0 \,.
\leqno(\frm{aa34})
$$
It is enough to show that
$$
\cp(\{\wtwoo>\delta\}\cap\{w=0\})=0
\leqno(\frm{aa34.5})
$$
for every $\delta>0$. If (\rf{aa34.5}) is not satisfied, by
Proposition~\rf{be01} we have ${\beta^{2}\over
\alpha^{2}} \mutwoo(\{{\wtwoo>\delta}\})=+\infty$, which
contradicts the fact that $\wtwoo\in L^{2}(\Om,\mutwoo)$. This
proves (\rf{aa34}).

As $w\ge0$ and $\woneo\ge 0$ q.e.\ on $\Om$ by the comparison principle
(Theorem 2.10
of [\rf{DM-Mos-86}]), from (\rf{aa33}) and (\rf{aa34}) it follows that
$$
\cp(\{\wtwoo>0\}\cap\{\woneo=0\})=
\cp(\{\wtwoo>0\}\setminus \{\woneo>0\})=0 \,,
\leqno(\frm{aa35})
$$
which proves (\rf{aa29++}).
Since $\nuoneo=\woneo\muoneo$ on $\{{\woneo>0}\}$ by (\rf{aa9}),
it follows from (\rf{aa35}) that  $\nuoneo=\woneo\muoneo$ on
$\{{\wtwoo>0}\}$, so that (\rf{aa31}) yields
$$
\woneo \muoneo \le {\beta^{2}\over \alpha^{2}}
\woneo \mutwoo \qquad \hbox{on }\, \{{\wtwoo>0}\}\,.
$$
As $\woneo>0$ q.e.\ on $\{{\wtwoo>0}\}$,  we conclude that
$$
\muoneo \le {\beta^{2}\over \alpha^{2}}
\mutwoo \qquad \hbox{on }\, \{{\wtwoo>0}\}\,.
\leqno(\frm{aa36})
$$

Let us finally prove that
$$
\muoneo \le {\beta^{2}\over \alpha^{2}}
\mutwoo \qquad \hbox{on }\, \{{\wtwoo=0}\}\,.
\leqno(\frm{aa37})
$$
Let $B$ be a Borel set contained in $\{{\wtwoo=0}\}$. If $\cp(B)=0$,
then $\muoneo(B)=\mutwoo(B)=0$, because $\muoneo$ and
$\mutwoo$ belong to $\mo$. If $\cp(B)>0$, then $\mutwoo(B)=+\infty$
by Proposition~\rf{be01}. In both cases we have
$\muoneo(B) \le {\beta^{2}\over \alpha^{2}}
\mutwoo(B)$, hence (\rf{aa37}) is proved.

Inequality (\rf{aa29})
now follows from (\rf{aa36}) and (\rf{aa37}).
\endproof

\parag{\chp{proofs}}{Proofs of the corrector results}

In this section we prove Lemma \rf{ab08.50} and Theorems \rf{ab061}
and \rf{ab07}, which give immediately
all results of Section \rf{corrector} (see the comments before the
statement of Theorem \rf{ab07}).

We begin by the following theorem, which is proved by
using the comparison result of Section~\rf{comparison}.

\th{\thm{af060}}{For every
$\e>0$, let $\ye\in H^{1}(\Om)\cap L^{2}(\Om,\mue)$. Assume that
$$
\ldisplaylinesno{
\ye\wto\yo \qquad\hbox{weakly in }\, \h\,,
&(\frm{af60})
\cr
\sup_{\e>0} \into |\ye|^{2} d\mue <+\infty \,.
&(\frm{af61})
\cr}
$$
Then $\yo\in L^{2}(\Om,\muo)$.
}

\proof We use the notion of $\gamma$-convergence, introduced in
[\rf{DM-Mos-87}] and further developed in [\rf{DM-87}],
which concerns the convergence of minima and minimizers of the
functionals $J^{\e}_{f}$ defined on $\xe$ by
$$
J^{\e}_{f}(y) = \alpha \int_{\Om} |Dy|^{2} dx +
\int_{\Om} |y|^{2} d\mue
- 2\langle f,y \rangle \,,
$$
for any given $f\in\hm$. Note that the minimizer of $J^{\e}_{f}$ is the
unique  solution to problem (\rf{bc1}) with $A=\alpha I$ and
$\mu=\mue$.
By Theorem~4.14 of [\rf{DM-Mos-87}] there exists a subsequence, still
denoted by $(\mue)$, which $\gamma$-converges (with respect to the
 operator $-\alpha\Delta$) to a measure
$\muoh\in\mo$ (the regularity property (b) of $\muoh$ is
obtained by using
Theorem~3.10 of [\rf{DM-87}]). By Lemma~5.5 of [\rf{DM-87}] we have
$$
\alpha \into |D\yo|^{2} dx + \into |\yo|^{2} d\muoh \le
\liminf_{\e\to0} \big\{
\alpha \into |D\ye|^{2} dx + \into |\ye|^{2} d\mue
\big\} < +\infty\,.
\leqno(\frm{af65})
$$

Let $\weh$ be the unique solution to
the problem
$$
\cases{\weh\in\xe\,,&
\cr\cr
\displaystyle
\alpha \int_{\Om} D\weh Dy\,dx +
\int_{\Om} \weh y \,d\mue
=\int_{\Om}y\,dx \qquad \forall y\in\xe\,.&
\cr}\leqno(\frm{af63})
$$
By Proposition~4.10 of [\rf{DM-Mos-87}] the sequence $(\weh)$
converges weakly in $\ho$ to the solution $\woh$ of the problem
$$
\cases{\woh\in\xoh\,,&
\cr\cr
\displaystyle
\alpha \int_{\Om} D\woh Dy\,dx +
\int_{\Om} \woh y \,d\muoh
=\int_{\Om}y\,dx \qquad \forall y\in\xoh\,.&
\cr}\leqno(\frm{af64})
$$
If we apply Theorem~\rf{aa04} with $\Aonee=\Ae$ and
$\Atwoe=\alpha I$, we obtain
$$
\muo\le {\beta^{2}\over \alpha^{2}}\muoh \,,
\leqno(\frm{ah14})
$$
so that (\rf{af65})
implies that $\yo\in L^{2}(\Om,\muo)$.
\endproof

\lemma{\thm{ab08.5}}{Under the hypotheses of Theorem~\rf{ab061},
we have
$$
\Ae D\ue \wto \Ao D\uo \qquad\hbox{weakly in }
\, \ltwon \,.
\leqno(\frm{ad22})
$$
Moreover there exists $\sigoo\in\mm\cap\hm$, with
$|\sigoo|\in\mm\cap\hm$, such that
$$
-\div(\Ao D\uo) + \sigoo=f \qquad \hbox{in }\, \distr\,.
\leqno(\frm{af2})
$$
}
\negbigskip

\proof Since the positive and the negative parts $(\ue)^{+}$ and
$(\ue)^{-}$ of $\ue$ belong to $\xe$, by Theorem
2.4 of [\rf{DM-Mos-86}] for every $\e>0$ we can consider
the solutions $\uonee$ and $\utwoe$ to the problems
$$
\ldisplaylinesno{
\hbox to 1cm{}
\cases{\uonee-(\ue)^{+} \in \xe\,,&
\cr\cr
\displaystyle
\into\Ae D\uonee Dy\,dx +
\into \uonee y \,d\mue
=\into f^{+} y\, dx
\qquad
\forall y\in \xe\,,&
\cr}&(\frm{ad50})
\cr\cr
\hbox to 1cm{}
\cases{\utwoe-(\ue)^{-} \in \xe\,,&
\cr\cr
\displaystyle
\into\Ae D\utwoe Dy\,dx +
\into \utwoe y \,d\mue
=\into f^{-} y\, dx
\qquad
\forall y\in \xe\,.&
\cr}&(\frm{ad51})
\cr}
$$
By linearity we have
$$
\ue=\uonee-\utwoe \quad\hbox{q.e.\ in }\,\Om\,.
\leqno(\frm{my1})
$$

Using $y=\uonee - (\ue)^{+}$ as test function in (\rf{ad50}), and then
(\rf{za12}) and (\rf{za13}),
as well as Poincar\'e's and Young's inequalities, we obtain
$$
\sup_{\e>0} \| \uonee \|_{\h} < +\infty \,.
\leqno(\frm{ad54})
$$
Passing to a
subsequence, we can assume that $(\uonee )$ converges
weakly in $\h$ to some function $\uoneo$.
Since $\uonee\ge0$ q.e.\ in $\Om$ by the comparison principle
(Theorem 2.10 of [\rf{DM-Mos-86}]), by
Proposition \rf{bc03} there exists
$\lape\in \hmp$ such that
$$
-\div(\Ae D\uonee) + \lape=f^{+}\quad \hbox{in }\, \distr\,.
\leqno(\frm{ad55})
$$
{}From Theorem~\rf{za006} we obtain that
$$
\Ae D\uonee \wto \Ao D\uoneo \qquad\hbox{weakly in }
\,L^{2}(\Om,\rn)\,,
\leqno(\frm{af4})
$$
and we deduce from (\rf{ad55}) that there exists
$\sigma^{0}_{\scriptscriptstyle\mskip-.5\thinmuskip\oplus}\in \hmp$ such
that
$$
-\div(\Ao D\uoneo) +
\sigma^{0}_{\scriptscriptstyle\mskip-.5\thinmuskip\oplus}
=f^{+}\qquad \hbox{in }\, \distr\,.
\leqno(\frm{af3})
$$

Properties (\rf{ad22}) and (\rf{af2}) now follow
from (\rf{af4}) and (\rf{af3}), from the analogous results for
$\utwoe$,  and  from (\rf{my1}).
\endproof

\proofof{Lemma~\rf{ab08.50}} By (\rf{my1}) we have $\ue=\uonee-\utwoe$
q.e.\ in $\Om$,
where $\uonee$ and $\utwoe$ are the solutions of (\rf{ad50}) and
(\rf{ad51}).
Let $\vonee$ be the solution to the problem
$$
\cases{\vonee-(\ue)^{+} \in \ho\,,&
\cr\cr
-\div(\Ae D\vonee)=f^{+} \quad\hbox{in }\, \distr\,.
&\cr}
$$
By the comparison principle (Theorem~2.10 of [\rf{DM-Mos-86}])
we have $0\le \uonee\le \vonee$ q.e.\ in $\Om$.

As $(\ue)^{+}$ is bounded in $\h$, the sequence $(\vonee)$ is bounded
in $\h$ too.
On the other hand the
classical local $L^{\infty}$ estimate for solutions of elliptic
equations
(see, e.g., [\rf{Sta}]) asserts that for every open set $U\subset\subset
\Om$
$$
\| \vonee \|_{L^{\infty}(U)} \le C_{U}
\| \vonee \|_{L^{2}(\Om)} \,,
\leqno(\frm{bi7})
$$
therefore $(\vonee)$, and hence $(\uonee)$, is bounded in
$L^{\infty}(U)$.
If $\ue\in\ho$, we have also $\vonee\in\ho$, and the
global $L^{\infty}$ estimate in $\Om$ implies that
$(\vonee)$, and hence $(\uonee)$, is bounded in $\linfty$.
A similar
argument holds for $(\utwoe )$, so that $(\ue)$ is bounded in
$L^{\infty}(U)$ (and also in $\linfty$
if $\ue\in\ho$) and
(\rf{bd1000}) is proved.

Let $\varphi$ be a function in $\cinftyo$ such that $\varphi=1$ in $U$.
Using
$y=\ue\varphi^{2}$ as test function in (\rf{ab23}), and then
(\rf{za12}), (\rf{za13}), and the boundedness of $(\ue)$ in $\h$,
we easily  obtain (\rf{bd1100}).
If  $\ue\in\ho$, we simply use $y=\ue$ as test function in (\rf{ab23}).
\endproof


\bigskip
The proof of Theorems \rf{ab061} and \rf{ab07} will be divided in
three lemmas.
For every $\e>0$ let $\ye$ be a function of $\h\cap\lmue$ such that
$$
\ye\wto\yo \qquad\hbox{weakly in }\h\,,
\leqno(\frm{ac1})
$$
for some function $\yo$ in $\h$. Assume that
$$
\ldisplaylinesno{
\sup_{\e>0} \|\ye\|_{\linfty} < +\infty \,,
&(\frm{ac2})
\cr
\sup_{\e>0} \into |\ye|^{2} d\mue <+\infty \,.
&(\frm{ac3})
\cr}
$$

\lemma{\thm{ad030}}{Under the hypotheses of Theorem~\rf{ab07},
let $\ye$, $\e\ge0$, be functions in $\h$ which satisfy
(\rf{ac1}), (\rf{ac2}), and (\rf{ac3}).
Then $\yo$ belongs to $\lmuo$ and for every $\varphi\in\cinftyo$ we
have
$$
\eqalign{
\lim_{\e\to0} &\into \Ae D\ve D\ye \varphi\,dx
+ \into \ve \ye\varphi \,d\mue=
\cr
&= \into \Ao D (\psi \,\omo) D\yo  \varphi\,dx  +
\into \yo \psi \, \omo \varphi \,d\muo \,.
\cr}\leqno(\frm{ac17})
$$
If, in addition, $\ye\in\ho$ for every $\e>0$, then (\rf{ac17}) also
holds
with $\varphi=1$.
}

\proof
We prove (\rf{ac17}) only in the case $\varphi\in\cinftyo$, since, under
the
additional hypothesis $\ye\in\ho$, the proof with $\varphi=1$ is
similar.
In this proof $(\se)$ will denote a sequence of real numbers converging
to $0$ as $\e$ tends to $0$, whose value can change from line to line.

Theorem~\rf{af060}, (\rf{ac1}), and (\rf{ac3}) imply that $\yo\in\lmuo$.

By (\rf{ab18.5}) for every $\varphi\in\cinftyo$ we have
$$
\eqalign{
\into & \Ae D\ve D\ye \varphi\,dx =
\cr
& =\sum_{j=1}^{n} \into \Dj{\psi} \Ae (\ej + D\zje ) D\ye \ome
\varphi\,dx
+ \into \psi\Ae D\ome D\ye \varphi\, dx + \se \,.
\cr}
$$
By (\rf{ba4}) and (\rf{ba5}) the sequence $(\ome)$ converges to
$\omo$ strongly in $L^{r}(\Om)$ for every $1\le r<+\infty$.
Since, by (\rf{ab5}), $(\zje)$ converges to $0$ weakly in in $\hp$
for some $p>2$,
we conclude that
$$
\into \Dj{\psi} \Ae (\ej + D\zje ) D\ye \ome
\varphi\,dx =
\into \Dj{\psi} \Ae (\ej + D\zje ) D\ye \omo
\varphi\,dx + \se\,.
$$
Therefore
$$
\into \Ae D\ve D\ye \varphi\,dx  + \into \ve \ye\varphi \,d\mue=
 I^{\e}+I\!I^{\e}+I\!I\!I^{\e}+\se \,,
\leqno(\frm{ac10})
$$
where
$$
\eqalign{
I^{\e}&=\sum_{j=1}^{n} \into \Dj{\psi} \Ae (\ej + D\zje ) D\ye \omo
\varphi\,dx \,,
\cr
I\!I^{\e}&= \into \psi\Ae D\ome D\ye \varphi\, dx \,,
\cr
I\!I\!I^{\e}&=\into \ve \ye\varphi \,d\mue \,.
\cr}
$$

We now pass to the limit in $ I^{\e}$, $I\!I^{\e}$, and
$I\!I\!I^{\e}$. For what concerns $ I^{\e}$, we write
$$
\eqalign{
 I^{\e}&=\into \Dj{\psi} \Ae (\ej + D\zje ) D\ye \omo \varphi\,dx=
\cr
&= \langle -\div( \Ae (\ej + D\zje )) , \Dj{\psi} \,\ye \omo \varphi
\rangle
- \into \Ae (\ej + D\zje ) D\Dj{\psi}\, \ye \omo \varphi \,dx -{}
\cr
&
- \into \Ae (\ej + D\zje ) \Dj{\psi} \,\ye D\omo \varphi
\,dx 
- \into \Ae (\ej + D\zje ) \Dj{\psi} \,\ye \omo D\varphi \,dx \,.
\cr}
$$
Properties (\rf{ab6}) and (\rf{ab2}) of $\zje$, together with
properties (\rf{ac1}) and (\rf{ac2}) of $\ye$, imply that we can pass
to the limit in each term of the right hand side of the previous
formula, so that
$$
 I^{\e}= \sum_{j=1}^{n} \into \Dj{\psi} \Ao \ej  D\yo \omo \varphi\,dx +
\se
=\into \Ao D\psi  D\yo \omo \varphi\,dx + \se\,.
\leqno(\frm{ac11})
$$

As for $I\!I^{\e}$, we write
$$
I\!I^{\e}= \into  \psi\Ae D\ome D\ye \varphi\, dx =
\into \Ae D\ome D(\ye \psi \,\varphi)\, dx
-\into \Ae D\ome\ye D(\psi \,\varphi)\, dx \,.
\leqno(\frm{ac13})
$$
As $\ome$ satisfies (\rf{ba1}), we have
$$
\into \Ae D\ome D(\ye \psi \,\varphi)\, dx +
\into \ome \ye \psi \,\varphi \,d\mue =
\into \ye \psi \,\varphi \,d\lae\,,
$$
and by (\rf{ba2.5}) and (\rf{ac1}) we conclude that
$$
\into \Ae D\ome D(\ye \psi \,\varphi)\, dx =
\into \yo \psi \,\varphi \,d\lao
-\into \ome \ye \psi \,\varphi \,d\mue + \se\,.
\leqno(\frm{ac14})
$$
By Proposition \rf{bc03} and Theorem \rf{za006}
$(\Ae D\ome)$ converges to $\Ao D\omo$ weakly in $\ltwon$,
while by (\rf{ac1}) $(\ye)$
converges to $\yo$ strongly in $\ltwoloc$.
Therefore
$$
-\into  \Ae D\ome\ye D(\psi \,\varphi)\, dx
=-\into  \Ao D\omo\yo D(\psi \,\varphi)\, dx + \se \,.
\leqno(\frm{ac15})
$$
{}From (\rf{ac13}), (\rf{ac14}), and (\rf{ac15}) we obtain that
$$
\eqalign{
I\!I^{\e}&= \into \yo \psi \,\varphi \,d\lao
-\into \ome \ye \psi \,\varphi \,d\mue -
\into  \Ao D\omo\yo D(\psi \,\varphi)\, dx + \se=
\cr
&= \into  \Ao D\omo D\yo \psi \,\varphi\, dx +
\into \omo \yo \psi \,\varphi \,d\muo
-\into \ome \ye \psi \,\varphi \,d\mue + \se  \,,
\cr}\leqno(\frm{ac16})
$$
where the last equality follows from (\rf{ba1}) for $\e=0$,
since $\yo\in\lmuo$.

Finally, we write $I\!I\!I^{\e}$ as
$$
I\!I\!I^{\e}=\into \ve\ye\varphi \,d\mue =
\into \psi\,\ome \ye\varphi \,d\mue +
\sum_{j=1}^{n} \into \Dj{\psi} \,\zje \ome \ye \varphi \,d\mue \,.
$$
Since, by (\rf{ab4}), $(\zje)$ converges to $0$ uniformly, while, by
(\rf{ba6}) and (\rf{ac3}), the norms of $\ome$ and $\ye$ in $\lmue$
remain bounded, we conclude that
$$
I\!I\!I^{\e} =
\into \psi\,\ome \ye\varphi \,d\mue + \se \,.
\leqno(\frm{ac12})
$$

{}From (\rf{ac10}), (\rf{ac11}), (\rf{ac16}), and (\rf{ac12})
we obtain (\rf{ac17}).
\endproof

\lemma{\thm{af04}}{Under the hypotheses of Theorem~\rf{ab07},
for every $\varphi\in\cinftyo$ we have
$$
\eqalign{
\lim_{\e\to0} &\big\{ \into \Ae D(\ue-\ve)D(\ue-\ve)\,\varphi\, dx
+\into |\ue-\ve|^{2}\varphi \,d\mue \big\} =
\cr
&= \into \Ao D(\uo-\psi\,\omo)D(\uo-\psi\,\omo)\,\varphi\, dx +{}
\cr
&\hphantom{=}+\into (\uo-\psi\,\omo) \,\varphi \,d\sigoo
-\into (\uo-\psi\,\omo) \psi\, \omo \varphi \,d\muo\,,
\cr}\leqno(\frm{af10})
$$
where $\sigoo$ is defined by (\rf{af2}).
If the functions $\ue$ and $\ve$ belong to $\ho$ for every $\e>0$,
then (\rf{af10}) also holds with $\varphi=1$.
}

\proof We prove the lemma only in the case $\varphi\in\cinftyo$,
since, under the additional hypothesis $\ue$, $\ve\in\ho$, the proof
with $\varphi=1$ is similar.

Let $\ye=\ue-\ve$ and let $\yo= \uo - \psi\omo$. Then
properties (\rf{ac1}), (\rf{ac2}), and
(\rf{ac3}) are satisfied by the
definition (\rf{af1}) of $\ve$ and by (\rf{ab3}), (\rf{ab4}),
(\rf{ba4}), (\rf{ba5}), (\rf{ba6}),
(\rf{ab20}), (\rf{bd10}), and (\rf{bd11}).
Using $y=\ye\varphi$ as test function in (\rf{ab23}) we
get
$$
\eqalign{
\into & \Ae D\ue D\ye \varphi\,dx
+ \into \ue \ye\varphi \,d\mue=
\cr
&= \into f \ye \varphi\,dx
- \into \Ae D\ue  D\varphi \,\ye\,dx \,.
\cr}\leqno(\frm{ad40})
$$
Using (\rf{ad22}) and (\rf{af2}) we obtain
$$
\eqalign{
\lim_{\e\to0} &\big\{ \into  \Ae D\ue D\ye \varphi\,dx
+ \into \ue \ye\varphi \,d\mue \big\} =
\cr
&= \into f \yo \varphi\,dx
- \into \Ao D\uo  D\varphi \,\yo\,dx =
\cr
&= \into \Ao D\uo  D\yo \varphi\,dx
+ \into \yo \varphi\,d\sigoo\,.
\cr}\leqno(\frm{af11})
$$
From (\rf{af11}) and (\rf{ac17}) we deduce (\rf{af10}).
\endproof

\lemma{\thm{af01}}{Under the hypotheses of Theorem~\rf{ab07},
for every
$y\in \xo$ we have
$$
\into y\,d\sigoo=\into y\,\uo \,d\muo\,,
\leqno(\frm{af30})
$$
where $\sigoo$ is defined by (\rf{af2}).
}

\proof First of all we recall that $\uo\in \Xo$ by Theorem \rf{af060}.
Let us fix $\varphi\in \cinftyo$ with $\varphi\ge 0$ in $\Om$.
By Lemma~\rf{af04} we have
$$
\eqalign{
\into & \Ao D(\uo-\psi\,\omo)D(\uo-\psi\,\omo)\,\varphi\, dx +{}
\cr
&{}+\into (\uo-\psi\,\omo) \,\varphi \,d\sigoo
-\into (\uo-\psi\,\omo)\, \psi\,\omo \varphi \,d\muo \ge 0\,,
\cr}\leqno(\frm{af41})
$$
for every $\psi\in\hh$.
By Proposition~\rf{bf02} the set
$\{\psi\,\omo:{\psi\in\cinftyo}\}$ is dense in $\xo$. Therefore
(\rf{af41}) implies that
$$
\into  \Ao D(\uo-z)D(\uo-z)\,\varphi\, dx +
\into (\uo-z)\, \varphi \,d\sigoo
-\into (\uo-z)\, z \,\varphi \,d\muo \ge 0\,,
\leqno(\frm{af42})
$$
for every $z\in\xo$.

We now use Minty's trick, and we take in (\rf{af42})
$z=\uo \zeta + t y$, with $t\in\r$, $y\in \xo$, and
$\zeta\in\cinftyo$ with $\zeta=1$ on $\supp\,\varphi$. Dividing by $t$
and
passing to the limit as $t$ tends to $0$ we obtain
$$
\into y \,\varphi \,d\sigoo =
\into y\, \uo \varphi \,d\muo\,.
$$
Since $\uo\in \lmuo$, we obtain (\rf{af30}) by approximating $1$ by a
sequence $(\varphi_k)$ of functions in
$\cinftyo$.
\endproof

\proofof{Theorem~\rf{ab061}} In view of Theorem \rf{af060}
the function $\uo$
belongs to $\h\cap\lmuo$. {}From (\rf{af2}) we have
$$
\into \Ao D\uo Dy \,dx + \into y\,d\sigoo
=\into f y \,dx
\qquad \forall y\in \ho\cap\lmuo  \,.
$$
By Lemma~\rf{af01} this implies (\rf{ab23}) for $\e=0$.
\endproof

\proofof{Theorem~\rf{ab07}} Since $\uo$
belongs to $\h\cap\lmuo$ by Theorem \rf{af060},
it is enough to apply Lemmas \rf{af04} and~\rf{af01}.
\endproof

\parag{\chp{Le}}{Problems with more general data}

In this section  we state and prove global and local convergence and
corrector results for relaxed Dirichlet problems of the form (\rf{bc30})
and
(\rf{ab23}), when the right hand sides $\fe$ and $f$ are replaced
by more general
linear functionals $\Le$, and when the strong convergence of $(\fe)$
in $\hm$ is replaced by the strong convergence of $(\Le)$ ``along the
sequence'' of spaces $(\xe)'$.

\goodbreak\bigskip
\noindent{\it Strong convergence of the data\/}\nobreak

For every $\e\ge0$ we consider an element of the dual space
$(\xe)'$, i.e., a linear functional $\Le\colon\xe\to\r$ such that
$$
|\Le(y)|\le C^{\e} \big\{ \alpha \into |Dy|^{2} dx + \into |y|^{2} d\mue
\big\}^{1\over2}
\qquad\forall y \in\xe \,,
$$
for a suitable constant $C^{\e}<+\infty$ (the constant $\alpha$ is
introduced in this formula for future convenience). It is easy to
prove that each functional $\Le$ can be represented in the form
$$
\Le(y)= \langle\fe,y\rangle + \into g^\e y\, d\mue\,,
\leqno(\frm{my11})
$$
where $\fe\in\hm$ and $g^\e\in\lmue$.

In this section we assume that
$$
\Le\to\Lo\qquad\hbox{strongly along the sequence }\, (\xe)'\,,
\leqno(\frm{my12})
$$
in the sense that
$$
\lim_{\e'\to0} L^{\e'}(y^{\e'}) = \Lo(\yo)\,,
\leqno(\frm{ah3})
$$
for every subsequence $\e'$ of $\e$ (see {\it Notation\/} in
Section~\rf{prel}) and every
sequence $(y^{\e'})$ which satisfies
$$
\ldisplaylinesno{
y^{\e'}\in \ho\cap L^{2}(\Om,\mu^{\e'}) \qquad \forall\e'>0\,,
&(\frm{ah4})
\cr
y^{\e'}\wto \yo \qquad \hbox{weakly in }\, \ho\,,
&(\frm{ah5})
\cr
\sup_{\e'>0} \into |y^{\e'}|^{2} d\mu^{\e'} <+\infty \,.
&(\frm{ah6})
\cr}
$$
Note that $\yo\in L^{2}(\Om,\muo)$ by Theorem~\rf{af060}.
Since (\rf{ah3}) holds true for every sequence $(y^{\e'})$ which
satisfies (\rf{ah4}),  (\rf{ah5}), and (\rf{ah6}), it is
easy to prove by contradiction that there exists a constant $C<+\infty$
such
that for every $\e>0$
$$
|\Le(y)|\le C \big\{ \alpha \into |Dy|^{2} dx + \into |y|^{2} d\mue
\big\}^{1\over2}
\qquad\forall y \in\xe \,.
\leqno(\frm{ah6.5})
$$

When $\Le$ is represented as in (\rf{my11}) with $g^\e=0$, it is easy
to see that (\rf{my12}) is satisfied if $(\fe)$ converges to $\fo$
strongly in $\hm$ (this condition is also necessary if all measures
$\mue$ are zero). The case where the functions $g^\e$ are not
identically zero is of course more difficult to handle, since the
measures $\mue$ vary, and the corresponding spaces $\lmue$ may be
different for different values of $\e$. This leads in a natural way
to definition (\rf{my12}), where we used the word ``strongly''
since the test functions $y^{\e'}$ in (\rf{ah3}) are only assumed to
be  uniformly bounded in
the corresponding spaces $\xe$.

In this definition the presence in (\rf{ah3}) of subsequences $\e'$
(and not just of the whole sequence $\e$) is due, among other
reasons, to the fact that we want that the convergence of $(\Le)$
implies the convergence of any subsequence.

\goodbreak\bigskip
\noindent{\it Global convergence and corrector results\/}\nobreak

By the Lax-Milgram lemma
for every $\e\ge0$ there exists a unique solution $\ue$
to the problem
$$
\cases{\ue\in\xe\,,&
\cr\cr
\displaystyle
\into\Ae D\ue Dy\,dx +
\into \ue y \,d\mue
= \Le(y) \qquad \forall y\in\xe\,.&
\cr}\leqno(\frm{ah7})
$$

The following theorem is a generalization of Theorem \rf{bg04}.

\th{\thm{ah02}}{Assume (\rf{ab1.5}), (\rf{ba1})--(\rf{bg50}), and
(\rf{my12}).
For every $\e\ge0$, let $\ue$ be the unique
solution to problem (\rf{ah7}).
Then $(\ue)$ converges to $\uo$ weakly in $\ho$.
}

\proof By (\rf{za12}), (\rf{za13}),  and (\rf{ah6.5}), using $y=\ue$ as
test
function in (\rf{ah7}) we obtain the estimate
$$
\alpha \into |D\ue|^{2} dx + \into |\ue|^{2} d\mue  \le  C^{2}\,.
\leqno(\frm{ah8.2})
$$
Extracting a subsequence, we may assume that
$$
\ue\wto \us \qquad \hbox{weakly in }\, \ho\,,
\leqno(\frm{ah8.4})
$$
for some function $\us\in\ho$.
By Theorem~\rf{af060} we have $\us\in L^{2}(\Om,\muo)$.
We will prove that $\us=\uo$. Since the limit does not depend on the
subsequence, this will prove that the whole sequence $(\ue)$ converges
to $\uo$.

If $y\in\xo$ satisfies $\into f y\,dx = 0$ for every
$f\in\linfty$, then $y=0$ a.e.\ in $\Om$. By the Hahn-Banach theorem,
this implies that $\linfty$ is dense in the dual space of $\xo$.
Therefore,
given $\deltaa>0$, there exists $\fd\in\linfty$ such that
$$
\big| \Lo(y)- \into \fd y\,dx \big| \le \deltaa
\big\{ \alpha \into |Dy|^{2} dx + \into |y|^{2} d\muo \big\}^{1\over2}
\quad\forall y \in\xo \,.
\leqno(\frm{ah9})
$$

For every $\e\ge0$ let $\ued$ be the unique solution to problem
(\rf{bc30}) with $\fe=\fd$.
By Theorem~\rf{bg04} we have
$$
\ued \wto \uod \qquad \hbox{weakly in }\, \ho\,,
\leqno(\frm{ah11})
$$
and taking $y=\ued$ as test function in (\rf{bc30}), with $\fe=\fd$, we
obtain
$$
\sup_{\e>0} \into |\ued |^{2} d\mue < +\infty \,.
\leqno(\frm{ah11.5})
$$

Using $y=\ue-\ued$ as test function in (\rf{ah7}) and (\rf{bc30}),
with $\fe=\fd$,  we obtain by difference
$$
\alpha \into |D(\ue-\ued)|^{2} dx + \into |\ue-\ued|^{2} d\mue
\le \Le(\ue-\ued) - \into \fd (\ue-\ued) \,dx \,,
\leqno(\frm{ah12})
$$
for every $\e\ge0$.
By (\rf{ah3}), (\rf{ah8.2}), (\rf{ah8.4}), (\rf{ah11}), and
(\rf{ah11.5}) we have
$$
\lim_{\e\to0} \big\{ \Le(\ue-\ued) + \into \fd (\ue-\ued) \,dx  \big\}
= \Lo(\us-\uod) - \into \fd (\us-\uod) \,dx  \,.
\leqno(\frm{ah13})
$$

Let $\muoh$ be the measure defined in the proof of Theorem~\rf{af060}.
By (\rf{af65}) we have
$$
\eqalign{
\alpha & \into |D(\us-\uod)|^{2} dx + \into |\us-\uod|^{2} d\muoh \le
\cr
& \le \liminf_{\e\to 0} \big\{
\alpha \into |D(\ue-\ued)|^{2} dx + \into |\ue-\ued|^{2} d\mue \big\}
\,.
\cr}\leqno(\frm{ah15})
$$
{}From (\rf{ah14}), (\rf{ah9}), (\rf{ah12}), (\rf{ah13}), and
(\rf{ah15}) we obtain that
$$
\alpha \into  |D(\us-\uod)|^{2} dx + \into |\us-\uod|^{2} d\muo \le
{\beta^4\over\alpha^4} \deltaa^{2}\,.
\leqno(\frm{ah16})
$$

Using (\rf{ah12}) for $\e=0$, we obtain from (\rf{ah9})
$$
\alpha \into  |D(\uo-\uod)|^{2} dx + \into |\uo-\uod|^{2} d\muo \le
\deltaa^{2}\,.
\leqno(\frm{ah17})
$$
{}From (\rf{ah16}) and (\rf{ah17}) we get
$$
\alpha \into  |D(\us-\uo)|^{2} dx\le
4 {\beta^4\over\alpha^4} \deltaa^{2}\,.
$$
Since $\deltaa>0$ is arbitrary, we conclude that $\us=\uo$.
\endproof

\bigskip

The next theorem is a generalization of Theorem \rf{ab03}.

\th{\thm{ah03}}{Assume (\rf{ab1.5}), (\rf{ba1})--(\rf{bg50}),
(\rf{ba5}), (\rf{ba6}), and (\rf{my12}).
Let $\delta>0$ and
let $\psid$ be a function in $\hh$ which satisfies (\rf{ab17}).
Assume that the functions $\ved$ defined by (\rf{ab18}) belong to
$\ho$. Then we have
$$
\limsup_{\e\to0} \big\{ \alpha \into |D\ue-D\ved|^{2} dx
+ \into |\ue-\ved|^{2} d\mue \big\} < \delta \,.
\leqno(\frm{ab1901})
$$
}
\negbigskip

\proof  Let us fix $\delta'<\delta$ such that
$$
\beta \into |D\uo-D(\psid\omo)|^{2} dx +
\into |\uo-\psid\omo|^{2} d\muo < \delta'\,.
\leqno(\frm{ab17+})
$$
For $\deltaa>0$, let $\fd$, $\ued$, and $\uod$ be as
in the proof of Theorem~\rf{ah02}.
Using (\rf{ab17+}) and (\rf{ah17}), we fix
$\deltaa>0$ small enough such that
$$
\ldisplaylinesno{
\sqrt{\delta'}+{\beta\over\alpha} \deltaa <\sqrt\delta \,,
&(\frm{ah20+})
\cr
\beta \into |D\uod-D(\psid\omo)|^{2} dx +
\into |\uod-\psid\omo|^{2} d\muo < \delta' \,.
&(\frm{ah20})
\cr}
$$
Therefore we can apply Theorem~\rf{ab03} with $f=\fd$ and we obtain
$$
\limsup_{\e\to0} \big\{ \alpha \into |D\ued-D\ved|^{2} dx
+ \into |\ued-\ved|^{2} d\mue \big\} < \delta'\,.
\leqno(\frm{ah21})
$$
As $(\ue)$ converges to $\uo$ weakly in $\ho$ by Theorem \rf{ah02},
using (\rf{ah9}), (\rf{ah12}), (\rf{ah13}), and (\rf{ah16}) we deduce
that
$$
\limsup_{\e\to0} \big\{ \alpha \into |D\ue - D\ued|^{2} dx
+ \into |\ue - \ued|^{2} d\mue \big\} \le {\beta^{2}\over \alpha^{2} }
\deltaa^{2} \,.
\leqno(\frm{ah22})
$$
{}From (\rf{ah20+}), (\rf{ah21}), and (\rf{ah22}) we obtain
(\rf{ab1901}).
\endproof

\goodbreak\bigskip
\noindent{\it Local convergence and corrector results\/}\nobreak

We consider now the case where the functions $\ue$ are solutions
to the problems
$$
\cases{\ue\in\Xe\,,&
\cr\cr
\displaystyle
\into\Ae D\ue Dy\,dx +
\into \ue y \,d\mue
= \Le(y) \qquad \forall y\in\xe\,,&
\cr}\leqno(\frm{ah701+})
$$
but are not required to satisfy the
boundary condition $\ue=0$ on $\partial\Om$.

The next theorem is a generalization of Corollary \rf{ab06}.

\th{\thm{ai01}}{Assume (\rf{ab1.5}), (\rf{ba1})--(\rf{bg50}),
(\rf{ba5}),
(\rf{ba6}), and (\rf{my12}). For every $\e>0$, let $\ue$ be a
solution to problem (\rf{ah701+}).
Assume that
$$
\ue\wto\uo \qquad\hbox{weakly in }\, \h\,,
\leqno(\frm{ab2001})
$$
for some function $\uo\in\h$.
Then $\uo$ is a solution to the problem
$$
\cases{
\uo\in\h\cap\lmuoloc\,,
&
\cr\cr
\displaystyle \int_{\Om}\Ao D\uo Dy\,dx +
\int_{\Om} \uo y \,d\muo
=\Lo(y) \qquad \forall y\in\hc\cap\lmuo\,,
&
\cr}
\leqno(\frm{ab2401})
$$
where $\hc$ denotes the space of all
functions $u\in \h$ with compact support in $\Om$.
If, in addition, $\uo\in\lmuo$, then the last line in (\rf{ab2401})
holds for every $y\in\xo$.
}

\proof As we have seen in the proof of Theorem \rf{ah02},
for every $\deltaa>0$ there exists $\fd\in\linfty$ which satisfies
(\rf{ah9}).
Let us fix an open set $U\subset\subset \Om$ and let
$\varphi\in\cinftyo$ such that $\varphi=1$ on $U$. Using
$y=\ue\varphi^2$ as test function in (\rf{ah701+}), and then
(\rf{za12}), (\rf{za13}), (\rf{ah6.5}), and (\rf{ab2001}) we obtain
$$
\sup_{\e>0} \intu |\ue |^{2} d\mue < +\infty \,,
\leqno(\frm{bd1111})
$$
which
implies that $\uo\in\lmuou$  by Theorem~\rf{af060}.
For every $\e\ge0$ let $\ued$ be the unique solution to the problem
$$
\cases{\ued-\ue\in\xeu\,,&
\cr\cr
\displaystyle
\intu\Ae D\ued Dy\,dx +
\intu \ued y \,d\mue
=\intu \fd y\,dx \qquad \forall y\in\xeu\,.&
\cr}\leqno(\frm{ah1001})
$$
Taking $y=\ued-\ue$ as test function in (\rf{ah1001}) and
(\rf{ah701+}),
we obtain by difference
$$
\alpha \intu |D(\ue-\ued)|^{2} dx + \intu |\ue-\ued|^{2} d\mue
\le \Le(\ue-\ued) - \intu \fd (\ue-\ued) \,dx \,,
\leqno(\frm{ah1201})
$$
for every $\e\ge0$.
By (\rf{ah6.5}) this implies that $(\ue-\ued)$ is bounded in $\hou$
and that the integrals $\intu |\ue-\ued|^{2} d\mue$ are bounded.
Using (\rf{ab2001}) and (\rf{bd1111}), we conclude that $(\ued)$
is bounded in $\hu$ and
$$
\sup_{\e>0} \intu |\ued|^{2} d\mue <+\infty \,.
\leqno(\frm{ai25})
$$
Extracting a subsequence, we may assume that
$$
\ued\wto \us \qquad \hbox{weakly in }\, \hu\,,
\leqno(\frm{ai30})
$$
for some function $\us\in\hu$ with $\us-\uo\in\hou$. By (\rf{ai25})
and by Theorem \rf{af060} the function $\us$ belongs to $\lmuou$.
Using both assertions of Corollary \rf{ab06}, $\us$ is a solution to the
problem
$$
\cases{
\us\in\hu\cap\lmuou\,,
&
\cr\cr
\displaystyle \intu\Ao D\us\, Dy\,dx +
\intu\us\, y \,d\muo
=\intu \fd y\,dx \qquad \forall y\in\hou\cap\lmuou\,.
&
\cr}
$$
Since $\us-\uo\in\xou$, by uniqueness, we have
$\us=\uod$.

By (\rf{my12}), (\rf{ab2001}), (\rf{bd1111}),
(\rf{ai25}), and (\rf{ai30}) we have
$$
\lim_{\e\to0} \big\{ \Le(\ue-\ued) - \intu \fd (\ue-\ued) \,dx  \big\}
= \Lo(\uo-\uod) - \intu \fd (\uo-\uod) \,dx  \,.
\leqno(\frm{ah1301})
$$

Let $\muoh$ be the measure defined in the proof of Theorem~\rf{af060}.
By (\rf{af65}) we have
$$
\eqalign{
\alpha & \intu |D(\uo-\uod)|^{2} dx + \intu |\uo-\uod|^{2} d\muoh \le
\cr
& \le \liminf_{\e\to 0} \big\{
\alpha \intu |D(\ue-\ued)|^{2} dx + \intu |\ue-\ued|^{2} d\mue
\big\} \,.
\cr}\leqno(\frm{ah1501})
$$
{}From (\rf{ah14}), (\rf{ah9}), (\rf{ah1201}), (\rf{ah1301}), and
(\rf{ah1501}) we obtain that
$$
\alpha \intu  |D(\uo-\uod)|^{2} dx + \intu |\uo-\uod|^{2} d\muo \le
{\beta^4\over\alpha^4} \deltaa^{2}\,.
\leqno(\frm{ah1601})
$$

Since,  by (\rf{ah9}), $\fd$ converges to $\Lo$ in the dual space
of $\xo$ as $\deltaa$ tends
to $0$, the solution $\uod$ of (\rf{ah1001}) for $\e=0$ converges in
$\xou$, as $\deltaa$ tends
to $0$, to the solution $v^{0}$ of the problem
$$
\cases{
v^{0}-\uo\in\xou\,,
&
\cr\cr
\displaystyle \intu \Ao Dv^{0} Dy\,dx +
\intu v^{0} y \,d\muo
=\Lo(y) \qquad \forall y\in\xou\,.
&
\cr}\leqno(\frm{ai40})
$$
On the other hand, by (\rf{ah1601}), $(\uod)$ converges to $\uo$ in
$\xou$
as $\deltaa$ tends to $0$. We conclude that $\uo=v^0$ and is the
solution
of
(\rf{ai40}). Since this holds  for every open set
$U\subset\subset\Om$, this implies that
$\uo$ is a solution of~(\rf{ab2401}).

The final statement of the theorem can be proved as explained before
Corollary~\rf{ab06}.
\endproof

\bigskip
The next theorem is a generalization of Theorem \rf{ad01}

\th{\thm{ad0101}}{Assume (\rf{ab1.5}), (\rf{ba1})--(\rf{bg50}),
(\rf{ba5}),
(\rf{ba6}), and (\rf{my12}). For every $\e>0$, let $\ue$ be a
solution to problem (\rf{ah701+}).
Assume that
$$
\ue\wto\uo \qquad\hbox{weakly in }\, \h\,,
$$
for some function $\uo\in\h$.
Let $U$ be an open set with
$U\subset\subset \Om$,
let $\delta>0$, let $\psid$ be a function in $\hhu$ which satisfies
(\rf{af51}), and let $\ved$ be the functions defined in $U$ by
(\rf{ab18}). Then
$$
\limsup_{\e\to0} \big\{ \alpha \intv  |D\ue-D\ved|^{2} dx
+ \intv  |\ue-\ved|^{2} d\mue \big\}< \delta \,,
\leqno(\frm{ad10+})
$$
for every open set $V\subset\subset U$.
}

\proof  Let us fix $\delta'<\delta$ such that
$$
\beta \intu |D\uo-D(\psid\omo)|^{2} dx +
\intu |\uo-\psid\omo|^{2} d\muo <  \delta' \,.
\leqno(\frm{af51+})
$$
For $\deltaa>0$, let $\fd$, $\ued$, and $\uod$ be as
in the proof of Theorem~\rf{ai01}.
Since $(\uod)$ converges to $\uo$ in $\xou$,
we fix $\deltaa$ small enough such that
$$
\ldisplaylinesno{
\sqrt{\delta'} + {\beta\over\alpha}\deltaa <\sqrt\delta
&(\frm{ah2001+})
\cr
\beta \intu |D\uod-D(\psid\omo)|^{2} dx +
\intu |\uod-\psid\omo|^{2} d\muo < \delta' \,.
&(\frm{ah2001})
\cr}
$$
Therefore we can apply Theorem~\rf{ad01} with $f=\fd$ and we obtain
$$
\limsup_{\e\to0} \big\{ \alpha \intv |D\ued-D\ved|^{2} dx
+ \intv |\ued-\ved|^{2} d\mue \big\} < \delta'\,,
\leqno(\frm{ah2101})
$$
for every open set $V\subset\subset U$.
Using (\rf{ah9}), (\rf{ah1201}), (\rf{ah1301}), and (\rf{ah1601})
we deduce that
$$
\limsup_{\e\to0} \big\{ \alpha \intu |D\ue - D\ued|^{2} dx
+ \intu |\ue - \ued|^{2} d\mue \big\} \le {\beta^{2}\over \alpha^{2} }
\deltaa^{2} \,.
\leqno(\frm{ah2201})
$$
{}From (\rf{ah2001+}), (\rf{ah2101}), and (\rf{ah2201}) we obtain
(\rf{ad10+}).
\endproof

\intro{References}
\ninepoint\frenchspacing

\item{[\bib{BLP}]}Bensoussan A., Lions J.L., Papanicolaou G.:
Asymptotic Analysis for Periodic Structures.
North Holland, Amsterdam, 1978.
\smallskip

\item{[\bib{BBDM}]}Bensoussan A., Boccardo L., Dall'Aglio A., Murat
F.: H-convergence for quasilinear elliptic equations under natural
hypotheses on the correctors.
{\it Proceedings of the Second Workshop on Composite Media and
Homogenization Theory (Trieste, 1993)\/}, {\it World Scientific,
Singapore\/}, 1995.
\smallskip

\item{[\bib{Boc-Mur 1}]}Boccardo L., Murat F.:
Nouveaux r\'esults de convergence dans des probl\`emes unilat\'eraux.
{\it Nonlinear Partial Differential Equations and Their Applications.
Coll\`ege de France Seminar. Vol. II\/}, 64-85 , {\it Res. Notes Math.,
Pitman,
London\/}, 1982
\smallskip

\item{[\bib{Boc-Mur 2}]}Boccardo L., Murat F.:
Increase of power leads to bilateral problems.
{\it Proceedings of the Second Workshop on Composite Media and
Homogenization Theory (Trieste, 1993)\/}, 113-123, {\it World
Scientific,
Singapore\/}, 1995.
\smallskip

\item{[\bib{But-DM-91}]}Buttazzo G., Dal Maso G.: Shape
optimization for Dirichlet problems: relaxed solutions and optimality
conditions. {\it Appl. Math. Optim.\/} {\bf 23} (1991), 17-49.
\smallskip

\item{[\bib{Cal-Cas}]}Calvo Jurado C., Casado Diaz J.:
The limit of Dirichlet systems for variable
monotone operators in general perforated domains.
{\it J. Math. Pures Appl.}, to appear.
\smallskip


\item{[\bib{Cio-Mur}]}Cioranescu D., Murat F.:
Un terme \'etrange venu d'ailleurs , I and II.
{\it Nonlinear Partial Differential Equations
and Their Applications. Coll\`ege de France Seminar. Vol. II\/},
98-138, {\it and Vol. III\/}, 154-178, {\it Res. Notes in Math., 60 and
70,
Pitman, London\/}, 1982 and 1983. English translation:
A strange term coming from nowhere.
{\it Topics in the Mathematical Modelling of Composite
Materials,\/} 45-93, {\it Birkh\"auser, Boston\/}, 1997.
\smallskip

\item{[\bib{Aglio-Mur}]}Dall'Aglio A., Murat F.:
A corrector result for $H$-converging parabolic problems with
time-dependent coefficients. {\it Ann. Scuola Norm. Sup. Pisa Cl. Sci.
(4)\/}
{\bf 25} (1997), 423-464.
\smallskip

\item{[\bib{DM-83}]}Dal Maso G.: On the integral representation of
certain local functionals. {\it Ricerche Mat.\/} {\bf 32} (1983),
85-113.

\item{[\bib{DM-87}]}Dal Maso G.:
$\Gamma$-convergence and  $\mu$-capacities.
{\it Ann. Scuola Norm. Sup. Pisa Cl. Sci. (4)\/}
{\bf 14} (1987), 423-464.
\smallskip

\item{[\bib{DM-Gar-94}]}Dal Maso G., Garroni A.:
New results on the asymptotic behaviour of Dirichlet problems
in perforated domains. {\it Math. Models Methods Appl. Sci.\/}
{\bf 4} (1994), 373-407.
\smallskip

\item{[\bib{DM-Mal}]}Dal Maso G., Malusa A.: Approximation of
relaxed Dirichlet problems by boundary value problems in perforated
domains. {\it Proc. Roy. Soc. Edinburgh Sect. A\/} {\bf 125} (1995),
99-114.

\item{[\bib{DM-Mos-86}]}Dal Maso G., Mosco U.:
Wiener criteria and energy decay for relaxed Dirichlet problems.
{\it Arch. Rational Mech. Anal.\/} {\bf 95} (1986), 345-387.
\smallskip

\item{[\bib{DM-Mos-87}]}Dal Maso G., Mosco U.:
Wiener's criterion and $\Gamma$-convergence.
{\it Appl. Math. Optim.\/} {\bf 15} (1987), 15-63.
\smallskip

\item{[\bib{DM-Mur-97}]}Dal Maso G., Murat F.:
Asymptotic behaviour and correctors for Dirichlet problems
in perforated domains with homogeneous monotone operators.
{\it Ann. Scuola Norm. Sup. Pisa Cl. Sci. (4)\/} {\bf 24} (1997),
239-290.
\smallskip

\item{[\bib{DM-Toa}]}Dal Maso G., Toader R.: Limits of Dirichlet
problems in perforated domains: a new formulation.
{\it Rend. Istit. Mat. Univ. Trieste\/} {\bf 24} (1994), 339-360.
\smallskip

\item{[\bib{Eva-Gar}]}Evans L.C., Gariepy R.F.: Measure Theory and
Fine Properties of Functions. CRC Press, Boca Raton, 1992.
\smallskip


\item{[\bib{Hei-Kil-Mar}]}Heinonen  J., Kilpel\"ainen T., Martio O.:
Nonlinear Potential Theory of Degenerate Elliptic Equations. Clarendon
Press, Oxford, 1993.
\smallskip

\item{[\bib{Kov}]}Kovalevsky A.:
An effect of double homogenization for Dirichlet problems in
variable domains of general structure.
{\it C. R. Acad. Sci. Paris S\'er. I Math.\/} {\bf 328} (1999)
1151-1156.
\smallskip

\item{[\bib{Maz}]}Maz'ya V.G.: Sobolev Spaces. Springer-Verlag,
Berlin, 1985.
\smallskip

\item{[\bib{Mey}]}Meyers N.G.:
An $L^p$-estimate for the gradient of
solutions of second order elliptic divergence equations.
{\it Ann. Scuola Norm. Sup. Pisa Cl. Sci. (3)\/}
{\bf 17} (1963), 189-206.
\smallskip

\item{[\bib{Mur-81}]}Murat F.:
L'injection du c\^one positif de $H^{?1}$ dans $W^{-1,q}$
est compacte pour tout $q < 2$. {\it J. Math. Pures Appl. (9)\/}
{\bf 60} (1981), 309-322.
\smallskip

\item{[\bib{Mur-Tar}]}Murat F., Tartar L.:
$H$-convergence. {\it S\'eminaire
d'Analyse Fonctionnelle et Num\'erique, Universit\'e d'Alger\/},
\hbox{1977-78}. English translation: Murat F., Tartar L.
$H$-convergence. {\it Topics in the Mathematical Modelling of
Composite Materials\/}, 21-43, {\it Birkh\"auser, Boston\/}, 1997.
\smallskip

\item{[\bib{Sanchez}]}Sanchez-Palencia E.:
Non Homogeneous Media and Vibration Theory.
Lecture Notes in Phys. 127, Springer-Verlag, Berlin, 1980.
\smallskip

\item{[\bib{Spa}]}Spagnolo S.:
Sulla convergenza di soluzioni di equazioni paraboliche ed ellittiche
{\it Ann. Scuola Norm. Sup. Pisa Cl. Sci. (3)\/} {\bf 22} (1968),
577-597.
\smallskip

\item{[\bib{Sta}]}Stampacchia G.: Le probl\`eme de Dirichlet pour
les \'equations elliptiques du second ordre \`a coefficients
discontinus.
{\it Ann. Inst. Fourier (Grenoble)\/} {\bf 15} (1965), 189-258.
\smallskip

\item{[\bib{Zie}]}Ziemer W.P.: Weakly Differentiable Functions.
Springer-Verlag, Berlin, 1989.
\smallskip

\bye